\documentclass[aap,MSNbibl,seceqn,dvips]{arximspdf}

\doi{10.1214/10-AAP741}
\volume{21}
\issue{5}
\pubyear{2011}
\firstpage{1827}
\lastpage{1859}

\makeatletter

\newtheorem{theorem}{Theorem}[subsection]
\newtheorem{lemma}[theorem]{Lemma}
\newtheorem{proposition}[theorem]{Proposition}
\newtheorem{corollary}[theorem]{Corollary}
\newproclaim{definition}[theorem]{Definition}
\newproclaim{example}[theorem]{Example}
\newproclaim{remark}[theorem]{Remark}
\newproclaim{remarks}{Remarks}

\makeatother

\begin{document}
\begin{frontmatter}

\title{Stochastic power law fluids:\\ Existence and uniqueness of
weak solutions}
\runtitle{Stochastic power law fluids}

\begin{aug}
\author[A]{\fnms{Yutaka} \snm{Terasawa}\ead[label=e1]{yutaka@ms.u-tokyo.ac.jp}}
\and
\author[B]{\fnms{Nobuo} \snm{Yoshida}\corref{}\thanksref{t1}\ead[label=e2]{nobuo@math.kyoto-u.ac.jp}
\ead[label=u1,url]{http://www.math.kyoto-u.ac.jp/\textasciitilde nobuo/}}
\thankstext{t1}{Supported in part by JSPS Grant-in-Aid for Scientific
Research, Kiban (C) 21540125.}
\runauthor{Y. Terasawa and N. Yoshida}
\affiliation{University of Tokyo and Kyoto University}
\address[A]{Graduate School of Mathematical Sciences\\
University of Tokyo\\
Komaba, Meguro-ku Tokyo 153-8914\\ Japan\\\printead{e1}} 
\address[B]{Division of Mathematics\\ Graduate School of Science\\
Kyoto University\\
Kyoto 606-8502\\ Japan\\
\printead{e2}\\\printead{u1}}
\end{aug}

\received{\smonth{2} \syear{2010}}
\revised{\smonth{7} \syear{2010}}

%
\begin{abstract}
We consider a stochastic partial differential equation (SPDE)
which describes the velocity field of
a viscous, incompressible non-Newtonian fluid
subject to a random force.
Here the extra stress tensor of the fluid is given by a polynomial of
degree $p-1$ of the rate of strain tensor, while
the colored noise is considered as a
random force. We investigate the existence and the
uniqueness of weak solutions to this SPDE.
\end{abstract}

%
\begin{keyword}[class=AMS]
\kwd[Primary ]{60H15}
\kwd[; secondary ]{76A05}
\kwd{76D05}.
\end{keyword}
\begin{keyword}
\kwd{Stochastic partial differential equation}
\kwd{power law fluids}.
\end{keyword}

\end{frontmatter}
%

\section{The power law fluids}\label{sec1}

We consider a viscous, incompressible fluid
whose motion is subject to a random force.
The container
of the fluid is supposed to be the torus
${\mathbb{T}}^d =({\mathbb{R}}/{\mathbb{Z}})^d \cong[0,1]^d$ as a
part of idealization.
For a differentiable vector field $v \dvtx{\mathbb{T}}^d \rightarrow
{\mathbb{R}}^d$, which is interpreted
as the velocity field of the fluid,
we denote the \textit{rate of strain tensor} by
\begin{equation}\label{e(v)}
e(v)
= \biggl( \frac{\partial_i v_j+\partial_j v_i }{2}\biggr)\dvtx
{\mathbb{T}}^d \rightarrow{\mathbb{R}}^d\otimes{\mathbb{R}}^d.
\end{equation}
We assume that
the extra stress tensor
\[
\tau(v) \dvtx {\mathbb{T}}^d \rightarrow{\mathbb{R}}^d\otimes{\mathbb{R}}^d
\]
depends on $e(v)$ polynomially. More precisely,
for $\nu>0$ (the kinematic viscosity) and $p>1$,
\begin{equation}\label{t(v)}
\tau( v ) = 2\nu\bigl(1 +|e(v)|^2\bigr)^{(p-2)/2}e(v).
\end{equation}
The linearly dependent case $p=2$ is the \textit{Newtonian fluid}
which is described by the Navier--Stokes equations, the special case
of (\ref{NS1}) and (\ref{NS2}).
On the other hand, both
the \textit{shear thinning} ($p<2$) and
the \textit{shear thickening} ($p>2$) cases are considered in many fields in
science and engineering. For example, shear thinning fluids are used
for automobile engine oil and pipeline for crude oil transportation,
while applications of shear thickening fluids can be found in
modeling of body armors and automobile four wheel driving systems.

Given an initial velocity $u_0\dvtx {\mathbb{T}}^d \rightarrow{\mathbb
{R}}^d$, the dynamics of the
fluid are described by
the following SPDE:
\begin{eqnarray}\label{NS1}
\operatorname{div}u&=&0,
 \\\label{NS2}
 \partial_tu+(u \cdot\nabla) u &=&-\nabla\Pi+
\operatorname{div}\tau(u) +\partial_t W,
\end{eqnarray}
where
\begin{equation}\label{divtau}
u \cdot\nabla=\sum_{j=1}^du_j \,\partial_j
    \quad\mbox{and}\quad
\operatorname{div}\tau(u)=\Biggl( \sum_{j=1}^d\partial_j\tau
_{ij}(u) \Biggr)_{i=1}^d.
\end{equation}
The unknown processes in the SPDE are the velocity field
$u =u(t,x)=(u_i(t,x))_{i=1}^d$ and the pressure
$\Pi=\Pi(t,x)$. The Brownian motion
$W =W(t,x)=(W_i(t,x))_{i=1}^d$ with values in $L_2 ({\mathbb{T}}^d
\rightarrow{\mathbb{R}}^d)$
(the set of vector fields on ${\mathbb{T}}^d$ with $L_2$ components)
is added as the random force.
Physical interpretations
of (\ref{NS1}) and (\ref{NS2}) are the mass conservation
and the motion equation, respectively. We note that the SPDE
(\ref{NS1}) and (\ref{NS2}) for the case
$p=2$ is the stochastic Navier--Stokes equation \cite{Fl08,FG95}.

Our motivation comes from works by
M\'alek et al.
\cite{MNRR96}, where the deterministic equation [the colored noise
$\partial_t W$ in (\ref{NS1}) and (\ref{NS2}) is replaced by a
nonrandom external force]
is investigated. Let
\begin{eqnarray}\label{p_1}
p_1(d)&=&\frac{3d }{ d+2} \vee\frac{3d -4 }{ d}=\cases{
\dfrac{3d }{ d+2}, &\quad \mbox{for $d \le4$},\vspace*{2pt}\cr
\dfrac{3d -4 }{ d}, &\quad \mbox{for $d \ge4$},
}
\\\label{p_2p_3}
p_2 (d)&=&\frac{2d }{ d-2},
   \qquad   p_3 (d)=\frac{3d-8 +\sqrt{9d^2+64} }{2d}
\end{eqnarray}
and
\begin{equation}\label{<p<}
p \in\cases{
( p_1(d), \infty),
&\quad \mbox{if $2 \le d \le8$},
\cr
(p_1(9),    p_2(9)) \cup(p_3(9),   \infty),
&\quad \mbox{if $d=9$},
\cr
(p_3 (d), \infty ), &\quad \mbox{if $d \ge10$}.
}
\end{equation}
For example, $p_1 (d)=\frac{3 }{2}$, $\frac{9 }{5}$,
$2$, $\frac{11 }{5}$ for $d=2,3,4,5$. A basic existence theorem
(\cite{MNRR96}, Theorem~3.4, page 222) states that the deterministic equation
has a weak solution if (\ref{<p<}) is satisfied, while a weak solution
is unique if $p \ge1+\frac{d }{2}$ (\cite{MNRR96}, Theorem~4.29, page 254).

The results in the present paper (Theorems~\ref{Thm.SNS} and~\ref{Thm.SNS1})
confirm that the above-mentioned
deterministic results are stable under
the random perturbation we consider.

Let us briefly sketch the outline of the proof of our existence result.

\begin{longlist}
\item[\textit{Step} $1$.] Set up a finite-dimensional subspace of a smooth,
divergence-free vector field, say $\mathcal{V}_n$, and an approximating
equation to the SPDE (\ref{NS1}) and (\ref{NS2}) in $\mathcal{V}_n$.
The good news here is that the approximating equation is a well posed
stochastic differential equation (SDE)
admitting a unique strong solution $u^n \in\mathcal{V}_n$. See
Theorem~\ref{Thm.Gar}
for detail.

\item[\textit{Step} $2$.] Establish some a priori bounds for the
solution $u^n \in\mathcal{V}_n$ of the approximating SDE
[e.g., (\ref{apriori}), (\ref{intb(X)a}), (\ref{intb(X)b}) and
(\ref{3.34})]. The point here
is that the bounds should be \textit{uniform in }$n$ for them to be
useful. Martingale
inequalities (e.g., the Burkholder--Davis--Gundy inequality) are
effectively used here, working in team with the Sobolev imbedding
theorem. See, for example, the proof of (\ref{apriori})  for details.

\item[\textit{Step} $3$.] Show that the solutions
$u^n \in\mathcal{V}_n$ to the approximating SDE are tight as $n
\rightarrow\infty$.
This is where the a priori bounds in step 2 play their roles as the moment
estimates to ensure that the tails of the solutions are thin enough
in certain Sobolev norms. This tightness argument is implemented
in Section \ref{conv_appro}.

\item[\textit{Step} $4$.] By step 3,
$u^n$ ($n \rightarrow\infty$) converges in law along a subsequence
to a limit. We verify that the limit is a weak solution to
the SPDE (\ref{NS1}) and (\ref{NS2}). These will be the subjects of Section
\ref{verify}.
\end{longlist}

Here are some comments concerning the technical difference
between the Navier--Stokes equations ($p=2$) and the power law fluids.
For the Navier--Stokes equations (both stochastic \cite{Fl08,FG95}
and deterministic \cite{Te79}),
it is reasonable to discuss solutions in the $L_2$-space.
On the other hand, for the power law fluids given by (\ref{t(v)}),
it is the $L_p$-space and its dual space that become relevant.
Also, due to the extra nonlinearity introduced by (\ref{t(v)}),
some of the arguments for $p \neq2$ become considerably
more involved than the case of $p=2$, especially for $p <2$.
(See, e.g., proof of Lemma~\ref{Lem.3.34}.)
We will overcome this difficulty by carrying the ideas in
\cite{MNRR96} over to the framework of It\^o's calculus.

\subsection{A weak formulation}

Let $\mathcal{V}$ be the set of ${\mathbb{R}}^d$-valued divergence
free, mean-zero
trigonometric polynomials, that is, the set of
$v \dvtx{\mathbb{T}}^d \rightarrow{\mathbb{R}}^d$ of the following
form:
\begin{equation}\label{cV}
v (x)=\sum_{z \in{\mathbb{Z}}^d\setminus\{ 0\}}\widehat{v}_z\psi
_z (x),   \qquad   x
\in{\mathbb{T}}^d,
\end{equation}
where
$\psi_z (x)=\exp(2 \pi\mathbf{i}z \cdot x)$ and the coefficients
$\widehat{v}_z \in{\mathbb{C}}^d$, $z \in{\mathbb{Z}}^d$ satisfy
\begin{eqnarray}\label{c_z =0}
\widehat{v}_z &= & 0\qquad
      \mbox{except for finitely many $z$,}
       \\ \label{v_real}
\overline{\widehat{v}_z} &= & \widehat{v}_{-z}
    \qquad  \mbox{for all $z$},
      \\\label{kc_z =0}
z \cdot\widehat{v}_z &= & 0
     \qquad \mbox{for all $z$}.
\end{eqnarray}
Note that (\ref{kc_z =0}) implies that
\[
\operatorname{div}v =0
    \qquad  \mbox{for all $v \in\mathcal{V}$.}
\]
For $\alpha\in{\mathbb{R}}$ and $v \in\mathcal{V}$ we define
\[
(1-\Delta)^{\alpha/2}v = \sum_{z \in{\mathbb{Z}}^d}(1+4\pi^2
|z|^2)^{\alpha/2}\widehat
{v}_z\psi_z.
\]
We equip the torus ${\mathbb{T}}^d$ with the Lebesgue measure.
For $p \in[1,\infty)$ and $\alpha\in{\mathbb{R}}$, we introduce
\begin{equation}\label{V_a}
V_{p,\alpha} = {}\mbox{the completion of $\mathcal{V}$ with respect to the
norm $ \|   \cdot  \|_{p,\alpha}$},
\end{equation}
where
\begin{equation}\label{normV_a}
\| v \|_{p,\alpha}^p=\int_{{\mathbb{T}}^d}|(1-\Delta)^{\alpha/2}v|^p.
\end{equation}
Then,
\begin{equation}\label{cpt_imb}
V_{p,\alpha+\beta} \subset V_{p,\alpha}
    \qquad  \mbox{for $1\le p <\infty$, $\alpha\in{\mathbb{R}}$ and
$\beta>0$}
\end{equation}
and
the inclusion $V_{p, \alpha+\beta} \rightarrow V_{p,\alpha}$ is compact
if $1 <p <\infty$ (\cite{Ta96}, (6.9), page 23).

For $v, w\dvtx{\mathbb{T}}^d \rightarrow{\mathbb{R}}^d$,
with $w$ supposed to be differentiable (for a moment),
we define a vector field
\begin{equation}\label{B(u,v)}
(v \cdot\nabla)w=\sum_{j}v_j\, \partial_j w
\end{equation}
which is bilinear in $(v,w)$. Later on, we will generalize the definition
of the above vector field; cf. (\ref{vwinV}).

Here are integration-by-parts formulae with which we reformulate
(\ref{NS1}) and (\ref{NS2}) into its weak formulation. In what follows,
the bracket $\langle u,v \rangle$ stands for the inner product of
$L_2 ({\mathbb{T}}^d \rightarrow{\mathbb{R}}^d)$, or its appropriate
generalization, for example,
the pairing of $u \in V_{p,\alpha}$ and $u \in V_{p',-\alpha}$
($p \in(1,\infty)$, $p'=\frac{p }{ p-1}$, $\alpha\ge0$). We let
$C^r ({\mathbb{T}}^d \rightarrow{\mathbb{R}}^d)$ ($r=1,\ldots,\infty$)
denote the
set of vector fields on ${\mathbb{T}}^d$ with $C^r$ components.

\begin{lemma}\label{Lem.uv=-vu}
For $v \in\mathcal{V}$ and $w,\varphi\in C^1 ({\mathbb{T}}^d
\rightarrow{\mathbb{R}}^d)$,
\begin{equation}\label{uv=-vu}
\langle\varphi, (v \cdot\nabla) w \rangle=-\langle w, (v \cdot
\nabla) \varphi\rangle.
\end{equation}
In particular,
\begin{equation}\label{uu=0}
\langle w, (v \cdot\nabla)w \rangle=0.
\end{equation}
Furthermore,
\begin{equation}\label{vp,t(v)}
\langle\varphi, \operatorname{div}\tau(v) \rangle=
- \langle\tau(v), e(\varphi) \rangle.
\end{equation}
\end{lemma}

\begin{pf} Since $\operatorname{div}v =0$, we have that
\[
\sum_{j}\partial_j (\varphi_i v_j)=
\sum_{j}\bigl( (\partial_j \varphi_i) v_j+\varphi_i \,\partial_j
v_j\bigr)
=\sum_{j}(\partial_j \varphi_i) v_j.
\]
Therefore,
\begin{eqnarray*}
\mbox{LHS of (\ref{uv=-vu})}
&=& \sum_{i,j}\langle\varphi_i, v_j\, \partial_j w_i \rangle
 =  -\sum_{i,j}\langle\partial_j (\varphi_i v_j), w_i \rangle
\\
& \stackrel{(1)}{=} &
-\sum_{i,j}\langle(\partial_j \varphi_i) v_j, w_i \rangle
=\mbox{RHS of (\ref{uv=-vu})}.
\end{eqnarray*}
Also, by integration by parts and the symmetry of $\tau_{ij}$,
\[
\mbox{LHS of (\ref{vp,t(v)})}=
-\sum_{i,j}\langle\partial_j \varphi_i, \tau_{ij} (v) \rangle
=-\sum_{i,j}\langle e_{ij} (\varphi), \tau_{ij} (v) \rangle
=\mbox{RHS of (\ref{vp,t(v)})}.\quad
\]
\upqed\end{pf}

Let us formally explain how the transformation of the
problem (\ref{NS1}) and (\ref{NS2}) into its weak formulation is achieved.
Suppose that $u,\Pi$ and ``$\partial_t W$'' in (\ref{NS1}) and (\ref{NS2})
are regular enough.
Then, for a test function $\varphi\in\mathcal{V}$,
%
\renewcommand{\theequation}{$\ast$}
\setcounter{equation}{0}
\begin{eqnarray}\label{*}
\quad&&\hspace*{20pt}\partial_t \langle\varphi, u \rangle=
-\underbrace{\langle\varphi, (u \cdot\nabla)u \rangle}_{(1)}
+ \underbrace{\langle\varphi, \operatorname{div}\tau(u)\rangle}_{(2)}-
\underbrace{ \langle\varphi, \nabla\Pi\rangle}_{(3)}+
\langle\partial_t W , \varphi\rangle,
\\
&&\hspace*{30pt}(1) \stackrel{\scriptsize{\mathrm{(\ref{uv=-vu})}}}{=}
-\langle(u \cdot\nabla)\varphi, u \rangle,
(2) \stackrel{\mbox{\scriptsize(\ref{vp,t(v)})}}{=}
- \langle e(\varphi), \tau(u) \rangle,\nonumber
(3) = -\langle\operatorname{div}\varphi,\Pi\rangle=0.\nonumber
\end{eqnarray}
Thus, ($\ast$) becomes
\[
\partial_t \langle\varphi, u \rangle=\langle(u \cdot\nabla
)\varphi, u \rangle
-\langle e(\varphi), \tau(u) \rangle+\partial_t \langle\varphi, W
\rangle.
\]
By integration, we arrive at
\renewcommand{\theequation}{\arabic{section}.\arabic{equation}}
\setcounter{equation}{19}
\begin{equation}\label{wNS}
\qquad\langle\varphi, u_t \rangle
=\langle\varphi, u_0 \rangle
+\int^t_0 \bigl(
\langle(u_s \cdot\nabla)\varphi, u_s \rangle-\langle e(\varphi
),\tau(u_s) \rangle
\bigr)\,ds+\langle\varphi, W_t \rangle.
\end{equation}
Here $u_t =u (t,   \cdot)$ and $W_t =W (t,   \cdot)$.
This is a standard weak formulation of (\ref{NS1}) and (\ref{NS2}).

\subsection{Bounds on the nonlinear terms}

Let us prepare a couple of $L_p$-bounds on the nonlinear terms.
They will be used to derive a priori bounds for the solutions later on.

\begin{lemma}\label{Lem.Buvw<}
Let $\alpha_i \in[0,\infty)$, $p_i \in[1,\infty )$, $i=1,2,3$, be
such that
\begin{equation}\label{sum_ia_i}
A \ge Bd,
   \qquad \mbox{where }A=\sum_i\alpha_i \mbox{ and } B=\sum_i\frac{1 }{ p_i}-1.\vspace*{-8pt}
\end{equation}
\begin{enumerate}[(a)]
\item[(a)]
Suppose $(\ref{sum_ia_i})$ and that $\frac{\alpha_i B }{ A}<\frac{1 }{ p_i}$
for all $i=1,2,3$.
Then, there exists
$C_1 \in(0,\infty)$ such that
\begin{equation}\label{Buvw<a}
| \langle w, (v \cdot\nabla)\varphi\rangle|
\le C_1\|v\|_{p_1,\alpha_1}\| w\|_{p_2,\alpha_2}\| \varphi\|
_{p_3,1+\alpha_3}
\end{equation}
for $v,w,\varphi\in C^\infty ({\mathbb{T}}^d \rightarrow{\mathbb{R}}^d)$.
\item[(b)] Suppose $(\ref{sum_ia_i})$, $\alpha_1+\alpha_2 >0$ and that
$B\le\frac{1 }{ p_i}$ for all $i=1,2,3$.
Then, for any $\theta\in(0,1)$, there exists
$C_2 \in(0,\infty)$ such that
\begin{equation}\label{Buvw<b}
| \langle w, (v \cdot\nabla)\varphi\rangle|
\le C_2
\|v\|_{p_1,\alpha_1}^\theta\|v\|_{p_1,\alpha_2}^{1 -\theta}
\| w\|_{p_2,\alpha_1}^{1 -\theta}\| w\|_{p_2,\alpha_2}^\theta\|
\varphi\|
_{p_3,1+\alpha_3}.
\end{equation}
\end{enumerate}
\end{lemma}

\begin{pf} (a)
Since
\[
\sum_{i,j}|w_i v_j\,\partial_j \varphi_i| \le|w||v||\nabla\varphi|,
\]
we have
\renewcommand{\theequation}{1}
\setcounter{equation}{0}
\begin{equation}
\quad | \langle w, (v \cdot\nabla)\varphi\rangle| \le
\|v \|_{q_1}\| w \|_{q_2}\| \nabla\varphi\|_{q_3}
\qquad \mbox{whenever } \frac{1 }{ q_1}+\frac{1 }{ q_2}+\frac{1 }{ q_3} \le1.
\end{equation}

\textit{Case} $1$. $B \le0$: We apply (1) with $q_i=p_i$ ($i=1,2,3$)
to get (\ref{Buvw<a}).

\textit{Case} $2$. $B>0$:
Since $\alpha\mapsto\| \cdot\|_{p_i,\alpha}$ is increasing
[$(1-\Delta)^{-\alpha/2}$ is a contraction on $L_p (T^d \rightarrow
{\mathbb{R}}^d)$
for any $\alpha\ge0$ and $p \ge1$],
it is enough to prove (\ref{Buvw<a}) with
$\alpha_i$ replaced by
$\tilde{\alpha}_i=\frac{\alpha_i }{ A}Bd.$
Therefore, we may assume without loss of generality that
\[
\max_i p_i \alpha_i <d
   \quad \mbox{and}  \quad  A = Bd.
\]
We apply (1) to $q_i \in[p_i,\infty)$, $i=1,2,3$ defined by
$\frac{1 }{ q_i}=\frac{1 }{ p_i}-\frac{\alpha_i }{ d}$.
We then use the following Sobolev imbedding theorem
(e.g., \cite{Ta96}, formula (2.11), page 5).
If~$\alpha p<d$ and $\frac{1 }{ q}=\frac{1 }{ p}-\frac{\alpha}{ d}$,
then there exists $C=C(d,\alpha) \in(0,\infty)$ such that
\renewcommand{\theequation}{\arabic{section}.\arabic{equation}}
\setcounter{equation}{23}
\begin{equation}\label{Sob}
\| v \|_q \le C \| v \|_{p,\alpha}\qquad
\mbox{for all $v \in C^\infty({\mathbb{T}}^d \rightarrow{\mathbb{R}}^d)$}.
\end{equation}

(b) Let us note the following interpolation inequality
(e.g., \cite{Ta96}, formula (6.5), page 23): for any $\lambda\in[0,1]$,
\renewcommand{\theequation}{2}
\setcounter{equation}{0}
\begin{equation}
\|u\|_{p_i,\lambda\alpha_1+(1-\lambda)\alpha_2 }
\le
C\|u\|_{p_i,\alpha_1}^\lambda\|u\|_{p_i,\alpha_2}^{1-\lambda}
\qquad \mbox{for } u \in V_{p_i,\alpha_1} \cap V_{p_i,\alpha_2}.
\end{equation}
On the other hand, we note that the assumptions for
(\ref{Buvw<a}) are satisfied if we replace $(\alpha_1,\alpha_2)$
by
\[
\bigl(\theta\alpha_1+(1-\theta)\alpha_2 , (1-\theta)\alpha_1+\theta
\alpha_2\bigr).
\]
Thus,
\begin{eqnarray*}
\qquad | \langle w, (v \cdot\nabla) \varphi\rangle|
&\stackrel{\mbox{\scriptsize(\ref{Buvw<a})}}{\le}&
C_1\|v\|_{p_1,\theta\alpha_1+(1-\theta)\alpha_2}
\|w\|_{p_2,(1-\theta)\alpha_1+\theta\alpha_2}\| \varphi\|
_{p_3,1+\alpha_3}\nonumber
\\
&\stackrel{\mbox{\scriptsize(2)}}{\le}& \mbox{RHS of (\ref{Buvw<b})}.\nonumber
\end{eqnarray*}
\upqed\end{pf}

\begin{lemma}\label{Lem.p2p}
Let $\alpha\in(0,1]$ and $p \in(\frac{2d }{ d+2\alpha}, \infty)$.
\begin{enumerate}[(a)]
\item[(a)]
Suppose that $(d,p,\alpha)\neq(2,2,1)$.
Then there exists
$C_1 \in(0,\infty)$ such that
\renewcommand{\theequation}{\arabic{section}.\arabic{equation}}
\setcounter{equation}{24}
\begin{equation}\label{Buvw<a2}
| \langle w, (v \cdot\nabla)\varphi\rangle|
\le C_1\|v\|_{p,\alpha}\| w\|_2\| \varphi\|_{p,\beta(p,\alpha)}
\end{equation}
for $v,w,\varphi\in C^\infty ({\mathbb{T}}^d \rightarrow{\mathbb
{R}}^d)$, where
\begin{equation}\label{a_4=}
\beta(p,\alpha) =\cases{
1+\biggl( \dfrac{2 }{ p}-\dfrac{1}{2}\biggr)d-\alpha>1,&\quad \mbox{if $ p < \dfrac{4d }{
d+2\alpha}$}, \cr
1, & \quad\mbox{if $p \ge\dfrac{4d }{ d+2\alpha}$}.
}
\end{equation}
\item[(b)]
Suppose that $d=2$. Then for any $\theta\in(0,1)$, there exists
$C_2 \in(0,\infty)$ such that
\begin{equation}\label{Buvw<b2}
| \langle w, (v \cdot\nabla)\varphi\rangle|
\le C_2
\|v\|_{2,1}^\theta\|v\|_2^{1 -\theta}
\| w\|_{p,1}^{1 -\theta}\| w\|_2^\theta\| \varphi\|_{2,1}
\end{equation}
for $v,w,\varphi\in C^\infty ({\mathbb{T}}^d \rightarrow{\mathbb{R}}^d)$.
\end{enumerate}
\end{lemma}

\begin{pf}
We apply Lemma~\ref{Lem.Buvw<} to
\[
(p_1,p_2,p_3)=(p,2,p), \qquad   (\alpha_1, \alpha_2)=(\alpha,0),\qquad
\alpha_3 = \biggl( \biggl(\dfrac{2 }{ p}-\dfrac{1 }{2}
\biggr)d-\alpha\biggr)^+.
\]
Then $\beta(p,\alpha)=1+\alpha_3$, $A=\alpha+\alpha_3$, $B=\frac{2
}{ p}-\frac{1 }{2}$.
It is enough to check that the assumptions of Lemma~\ref{Lem.Buvw<}(b) are
satisfied if $(d,p,\alpha) = (2,2,1)$ and that the assumptions of
Lemma~\ref{Lem.Buvw<}(a)
are satisfied if $(d,p,\alpha) \neq(2,2,1)$. In fact, the verification
for the case $(d,p,\alpha) = (2,2,1)$ can be done by simply plugging
the values.
We assume $(d,p,\alpha) \neq(2,2,1)$ in what follows.
We may assume that \mbox{$B>0$}, or equivalently $p<4$.
We have $A \ge Bd$ by the choice of $\alpha_i$'s.
Let us check that
\renewcommand{\theequation}{1}
\setcounter{equation}{0}
\begin{equation}
\frac{\alpha_1 }{ A}B=\frac{\alpha}{\alpha+\alpha_3}
\biggl(\frac{2 }{ p}-\frac{1 }{2}\biggr)<\frac{1 }{ p}.
\end{equation}
If $(d,p,\alpha) \neq(2,2,1)$ and
$p \ge\frac{4d }{ d+2\alpha}$ (which implies $p>2$),
then $\alpha_3=0$ and (1) is satisfied. If $(d,p,\alpha) \neq
(2,2,1)$ and
$p < \frac{4d }{ d+2\alpha}$ (which implies $p< \frac{d }{\alpha}$),
then $\alpha_3=( \frac{2 }{ p}-\frac{1 }{2}
)d-\alpha>0$.
One then sees that (1) is equivalent to that $p<\frac{d }{\alpha}$
and hence, is satisfied.
Let us check that
\renewcommand{\theequation}{2}
\setcounter{equation}{0}
\begin{equation}
\frac{\alpha_3 }{ A}B=\frac{\alpha_3 }{\alpha+\alpha_3}
\biggl(\frac{2 }{ p}-\frac{1 }{2}\biggr)<\frac{1 }{ p}.
\end{equation}
If $(d,p,\alpha) \neq(2,2,1)$ and
$p \ge\frac{4d }{ d+2\alpha}$,
then $\alpha_3=0$ and (2) is satisfied. If $p < \frac{4d }{ d+2\alpha}$,
then $\alpha_3=( \frac{2 }{ p}-\frac{1 }{2}
)d-\alpha>0$.
One then sees that (2)
is equivalent to that $p >\frac{2d }{ d+2\alpha}$ and
hence, is satisfied.
\end{pf}

\begin{remark*}
We note that the following variant of (\ref{Buvw<a2})
is also true:
\renewcommand{\theequation}{\arabic{section}.\arabic{equation}}
\setcounter{equation}{27}
\begin{equation}\label{Buvw<a3}
| \langle w, (v \cdot\nabla)\varphi\rangle|
\le C_1\|v\|_2\| w\|_{p,\alpha}\| \varphi\|_{p,\beta(p,\alpha)}.
\end{equation}
This can be seen by
interchanging the role of $(p_1,\alpha_1)$ and $(p_2,\alpha_2)$
in the above proof.
\end{remark*}

\begin{lemma}\label{Lem.b_2<}
For $p \in(1,\infty)$, there exists
$C_1 \in(0,\infty)$ such that
\begin{eqnarray}\label{b_2<}
&&| \langle e(\varphi), \tau(v) \rangle|
\le C_1\bigl(1+\|e (v)\|_p\bigr)^{p-1}\| e(\varphi) \|_p\nonumber
   \\[-8pt]\\[-8pt]
    &&\eqntext{\mbox{for all $v \in V_{p,1}$ and
    $\varphi\in\mathcal{V}$.}}
\end{eqnarray}
\end{lemma}

\begin{pf} Since
\[
|\tau(v)| \le C \bigl(1+|e (v)|\bigr)^{p-1},
\]
we have that
\begin{eqnarray*}
| \langle e(\varphi), \tau(v) \rangle|
&\le& C \int_{{\mathbb{T}}^d}\bigl(1+|e (v)|\bigr)^{p-1}|e(\varphi)|
\\
& \stackrel{{(p-1) / p}+{1 / p}=1}{\le} &
C \bigl\|1+|e (v)| \bigr\|_p^{p-1}\| e(\varphi) \|_p
\\
& \le& C \bigl(1+\|e (v)\|_p\bigr)^{p-1}\| e(\varphi) \|_p,
\end{eqnarray*}
which proves (\ref{b_2<}).
\end{pf}

Let $p \in(\frac{2d }{ d+2}, \infty)$, $v,w \in V_{p,1} \cap V_{2,0}$
and $u \in V_{p,1}$.
In view of Lemma~\ref{Lem.uv=-vu}, we think of
$(v \cdot\nabla)w$ and
$ \operatorname{div}\tau(u) $, respectively, as the following linear
functionals on $\mathcal{V}$:
\begin{eqnarray*}
 \varphi&\mapsto&\langle\varphi, (v \cdot\nabla)w\rangle
\stackrel{\rm def.}{=}-\langle w, (v \cdot\nabla)\varphi\rangle,
\\
 \varphi&\mapsto&\langle\varphi, \operatorname{div}\tau(u)
\rangle
\stackrel{\rm def.}{=}
- \langle e(\varphi), \tau(u) \rangle.
\end{eqnarray*}
Then, by Lemmas~\ref{Lem.p2p} and \ref{Lem.b_2<}, they
extend continuously, respectively, on $V_{p,\beta(p,1)}$ and
on $V_{p,1}$, where
\begin{equation}\label{a_4=1}
\beta(p,1) =\cases{
\biggl( \dfrac{2 }{ p}-\dfrac{1}{2}\biggr)d >1,&\quad \mbox{if $ p < \dfrac{4d }{ d+2}$},
\cr
1, &\quad \mbox{if $p \ge\dfrac{4d }{ d+2}$};
}
\end{equation}
cf. (\ref{a_4=}).
This way, we regard $(v \cdot\nabla) w \in V_{p',-\beta(p,1)}$
($p'=\frac{p }{ p-1}$)
with
\begin{eqnarray}\label{vwinV}
&&\| (v \cdot\nabla)w \|_{p',-\beta(p,1)}\nonumber
\\[-8pt]\\[-8pt]
&&\qquad \le
\cases{
C\|v\|_{2,1}^\theta\|v\|_{2}^{1 -\theta}
\|w\|_{2,1}^{1 -\theta}\| w\|_{2}^\theta, & \quad\mbox{if $p=d=2$},\cr
C \|v\|_{p,1}\| w\|_{2},
&\quad\mbox{if otherwise},
}\nonumber
\end{eqnarray}
and $\operatorname{div}\tau(u) \in V_{p',-1}$ with
\begin{equation}\label{6tauinV}
\| \operatorname{div}\tau(u) \|_{p',-1} \le
C \bigl(1+\|e(u)\|_p\bigr)^{p-1}.
\end{equation}
Finally, for $v \in V_{p,1} \cap V_{2,0}$, we define
\begin{equation}\label{b(v)}
b(v)=-(v \cdot\nabla)v+\operatorname{div}\tau(v) \in V_{p',-\beta
(p,1)}.
\end{equation}
With this notation, (\ref{wNS}) takes the form
\[
\langle\varphi, u_t \rangle
=\langle\varphi, u_0 \rangle
+\int^t_0 \langle\varphi, b(u_s)\rangle\, ds
+\langle\varphi, W_t \rangle,
\]
that is,
\begin{equation}\label{wNS2}
u_t =u_0+\int^t_0 b(u_s)\,ds+W_t
\end{equation}
as linear functionals on $\mathcal{V}$.

\section{The stochastic power law fluids} \label{Splf}

\subsection{The existence theorem}

We need the following definition.

\begin{definition}\label{Def.HBM}
Let $H$ be a Hilbert space and
$\Gamma\dvtx H \rightarrow H$ be a self-adjoint, nonnegative definite
operator of trace class.
A random variable $(W_t)_{t \ge0}$ with values in $C([0,\infty)
\rightarrow H)$
is called an {\it$H$-valued Brownian motion} with the covariance
operator $\Gamma$ [abbreviated by BM($H,\Gamma$) below]
if, for each $\varphi\in H$ and $0 \le s <t$,
\[
E [ \exp(\mathbf{i}\langle\varphi, W_t -W_s\rangle
)| (W_u)_{u \le
s} ]
=\exp\biggl(-\frac{t-s }{2}\langle\varphi, \Gamma\varphi\rangle
\biggr),
  \qquad  \mbox{a.s.}
\]
\end{definition}

To introduce the notion of weak solution (Definition~\ref{Def.w_sol}),
we agree on the following standard notation and convention.
For a Banach space $X$, we let
$L_{q,{\rm loc}} ([0,\infty) \rightarrow X)$ ($1 \le q \le\infty$)
denote the
set of locally $L_q$-functions $u\dvtx [0,\infty) \rightarrow X$,
with the Fr\'echet space metric induced by the semi-norms
$\| u \|_{L_q ([0,T] \rightarrow X)}$, $0<T <\infty$, where
$\| u \|_{L_q ([0,T] \rightarrow X)}$ stands for the standard
$L_q$-norm for
$u|_{[0,T]}\dvtx [0,T] \rightarrow X$. We also regard $C ([0,\infty)
\rightarrow X)$,
the set of continuous functions $u\dvtx [0,\infty) \rightarrow X$,
as the Fr\'echet space
induced by the semi-norms
$\sup_{0 \le t \le T}\| u(t) \|_{X}$, $0<T <\infty$.

We recall that the number $p$ is from (\ref{t(v)})
and that $b(v) \in V_{p', -\beta(p,1)}$ for
$v \in V_{p,1} \cap V_{2,0}$ is defined by (\ref{b(v)}).

\begin{definition}\label{Def.w_sol}
Suppose that:
\begin{itemize}[$\bullet$]
\item[$\bullet$]
$\Gamma\dvtx V_{2,0} \rightarrow V_{2,0}$ is a bounded self-adjoint,
nonnegative definite
operator of trace class;

\item[$\bullet$]
$\mu_0$ is a Borel probability measure on $V_{2,0}$;

\item[$\bullet$]
$(X,Y)=( (X_t,Y_t))_{t \ge0}$ is a process defined on a
probability space $(\Omega, \mathcal{F}, P)$ such that
\begin{eqnarray}\label{state}
X &\in& L_{p,{\rm loc}} \bigl([0,\infty) \rightarrow V_{p,1}\bigr)
\cap L_{\infty,{\rm loc}} \bigl([0,\infty) \rightarrow V_{2,0}\bigr)\nonumber
\nonumber\\[-8pt]\\[-8pt]
&&{}\cap C \bigl([0,\infty) \rightarrow V_{2\wedge p',-\beta}\bigr)\nonumber
\end{eqnarray}
for some $\beta>0$ and
$(Y_t)_{t \ge0}$ is a ${\rm BM}(V_{2,0}, \Gamma)$; cf.
Definition~\ref{Def.HBM}.
\end{itemize}

 Then the process $(X,Y)$ is said to
be a \textit{weak solution} to the SDE
\begin{equation}\label{SNS_def}
X_t = X_0+\int^t_0 b(X_s)\,ds +Y_t
\end{equation}
with the initial law $\mu_0$ if the following conditions
are satisfied:
\begin{eqnarray}\label{X_0=m}
& P(X_0 \in\cdot)=\mu_0 ;  &
\\\label{X_0indepM}
& \quad
 Y_{t+\cdot}-Y_t \quad\mbox{and}\quad
\{ \langle\varphi, X_s \rangle  ;   s \le t, \varphi\in\mathcal
{V}\}
\quad\mbox{are independent for any $t \ge0$};\hspace*{-25pt}&
\\\label{SNS}
& \quad
\langle\varphi, X_t \rangle= \langle\varphi, X_0 \rangle+\displaystyle\int^t_0
\langle\varphi, b(X_s) \rangle\, ds +\langle\varphi, Y_t \rangle\quad
\mbox{for all $\varphi\in\mathcal{V}$ and $t \ge0$}.\hspace*{-27pt}&
\end{eqnarray}
\end{definition}

We can now state our existence result.

\begin{theorem}\label{Thm.SNS}
Let $\Gamma$ and $\mu_0$ be as in Definition~$\ref{Def.w_sol}$ and
suppose additionally that:
\begin{itemize}[$\bullet$]
\item[$\bullet$] $(\ref{<p<})$ holds;

\item[$\bullet$]
$\Delta\Gamma=\Gamma\Delta$ and both $\Gamma$, $\Delta\Gamma$ are
of trace class;

\item[$\bullet$]
$\mu_0$ is a probability measure on $V_{2,1}$ and
\begin{equation}\label{m_alpha}
m_\alpha=\int\| \xi\|_{2,\alpha}^2\mu_0 (d\xi)<\infty  \qquad  \mbox
{for $\alpha=0,1$}.
\end{equation}
\end{itemize}
 Then there exists a weak solution
to the SDE $(\ref{SNS_def})$
with the initial law $\mu_0$;  cf. \textit{Definition}~\ref{Def.w_sol}
 such that
$(\ref{state})$ holds with $\beta=\beta(p,1)$;  cf. $(\ref{a_4=1})$.
Moreover, for any $T>0$,
\begin{equation}\label{apriori2}
E \biggl[ \sup_{t \le T}\| X_t \|_2^2 +\int^T_0\| X_t \|_{p,1}^p\, dt
\biggr]
\le(1+T)C <\infty,
\end{equation}
where $C=C(d,p,\Gamma,m_0) <\infty$.
\end{theorem}

\begin{remark*}
It would be worthwhile to mention that
Theorem~\ref{Thm.SNS} with $p=2$ is valid for \textit{all} $d$,
although it is not covered by the condition (\ref{<p<})
if $d \ge4$. In fact, Lemma~\ref{Lem.3.34}
is the only place we need condition (\ref{<p<}).
For $p=2$, however, we can avoid the use
of that lemma; cf. remarks at the end of Section
\ref{conv_appro} and after Lemma~\ref{Lem.conv_intb}.
\end{remark*}

\subsection{The uniqueness theorem}

As in the case of the deterministic equation \cite{MNRR96}, Theorem 4.29,
page 254,
we have the following uniqueness result:

\begin{theorem}\label{Thm.SNS1}
Suppose that
\begin{equation}\label{p>uni}
p \ge1+\frac{d }{2}.
\end{equation}
Then the weak solution
to the SDE $(\ref{SNS_def})$, subject to the  {a priori}  bound
$(\ref{apriori2}),$ is pathwise unique in the
following sense: if $(X,Y)$ and $(\widetilde{X},Y)$ are two solutions
on a
common probability space $(\Omega, \mathcal{F}, P)$ with a
common $\operatorname{BM}(V_{2,0}, \Gamma)$ $Y$ such that $X_0=\widetilde{X}_0$ a.s.,
then,
\[
P (X_t =\widetilde{X}_t  \mbox{ for all $t \ge0$})=1.
\]
%
\end{theorem}

The above uniqueness theorem, together with the Yamada--Watanabe
theorem provides us with
the so-called \textit{strong solution in the stochastic sense} to the SDE~(\ref{SNS_def}).

\begin{corollary}\label{Cor.YW}
Suppose $(\ref{p>uni})$, in addition to all the assumptions in
$Theorem~\ref{Thm.SNS}$, and let
$\xi$ be a given $V_{2,0}$-valued random variable with the law $\mu_0$
and $Y$ be a
given $\operatorname{BM}(V_{2,0}, \Gamma)$ independent of $\xi$.
Then there exists a process $X$ obtained as a function of $(\xi,Y)$,
such that $(X,Y)$ is weak solution
to the SDE $(\ref{SNS_def})$ with $X_0=\xi$ and with
all the properties stated in $Theorem~\ref{Thm.SNS}$. Moreover, the
law of the above
process $X$ is unique.
\end{corollary}

\begin{pf} Corollary~\ref{Cor.YW} is a direct consequence of Theorems~\ref{Thm.SNS}
and \ref{Thm.SNS1} via the Yamada--Watanabe theorem
\cite{IW89}, Theorem 1.1, page~163. The Yamada--Watanabe theorem
is usually stated for SDEs in finite dimensions.
However, as is obvious from its proof, it applies to the present setting.
\end{pf}

\begin{remark}
For $p \in[1+\frac{d }{2}, \frac{2d }{ d-2})$,
an even stronger version of Corollary~\ref{Cor.YW} is shown in \cite{Yo10}
as a consequence of strong convergence of the Galerkin approximation;
cf. Section \ref{sec:Gal}.
\end{remark}

\section{The Galerkin approximation} \label{sec:Gal}

\subsection{The exsitence theorem for the approximations}

For each $z \in{\mathbb{Z}}^d\setminus\{ 0\}$, let
$\{e_{z,j}\}^{d-1}_{j=1}$
be an orthonormal basis of the hyperplane
$\{ x \in{\mathbb{R}}^d;   z \cdot x =0\}$ and let
\begin{eqnarray}\label{psi}
 \quad&&\psi_{z,j}(x)\nonumber
\\[-8pt]\\[-8pt]
 &&\qquad=
\cases{
\sqrt{2}e_{z,j}\cos(2\pi z \cdot x), &\quad $j=1,\ldots,d-1$, \cr
\sqrt{2}e_{z,j -d+1}\sin(2\pi z \cdot x), &\quad $j=d,\ldots,2d-2$,
}
\qquad     x \in{\mathbb{T}}^d.\nonumber
\end{eqnarray}
Then
\[
\bigl\{ \psi_{z,j};     (z,j) \in({\mathbb{Z}}^d\setminus\{ 0\})
\times\{
1,\ldots,2d-2\} \bigr\}
\]
is an orthonormal basis of $V_{2,0}$.
We also introduce
\begin{eqnarray}\label{cV^n}
\mathcal{V}_n &= & {}\mbox{the linear span of }
\{\psi_{z,j};   (z,j) \mbox{ with } z \in[-n,n]^d\};\nonumber
\\[-8pt]\\[-8pt]
\mathcal{P}_n & = & \mbox{the orthogonal projection}{}\dvtx V_{2,0}
\rightarrow\mathcal{V}_n.\nonumber
\end{eqnarray}
Using the orthonormal basis (\ref{psi}),
we identify $\mathcal{V}_n$ with ${\mathbb{R}}^N$, $N=\dim\mathcal{V}_n$.
Let $\mu_0$ and $\Gamma, V_{2,0} \rightarrow V_{2,0}$, be as in
Theorem~\ref{Thm.SNS}.
Let also $\xi$ be a random variable such that
$P (\xi\in\cdot)=\mu_0$.
Finally, let $W_t$ be a BM$(V_{2,0},\Gamma)$ defined on a probability
space $(\Omega^W, \mathcal{F}^W, P^W)$. Then, $\mathcal{P}_n W_t$ is
identified with
an $N$-dimensional Brownian motion with covariance matrix $\Gamma
\mathcal{P}_n$.
Then we consider the following
approximation of (\ref{SNS}):
\begin{equation}\label{sde_Gar^n}
X_t^n=X_0^n+\int^t_0 \mathcal{P}_n b (X^n_s)\,ds+\mathcal{P}_n W_t,
   \qquad   t \ge0,
\end{equation}
where $X_0^n=\mathcal{P}_n\xi$. Let
\begin{equation}\label{X^zj}
X^{n,z,j}_t=\langle X^n_t, \psi_{z,j} \rangle
\end{equation}
be the $(z,j)$-coordinate of $X^n_t$.
Then (\ref{sde_Gar^n}) reads
\begin{equation}\label{sde_Gar^n2}
X^{n,z,j}_t=X^{n,z,j}_0 +\int^t_0 b^{z,j} (X^n_s)\,ds
+W^{z,j}_t,
\end{equation}
where
\begin{eqnarray}\label{b^zj}
b^{z,j} (X^n_s)&=&\langle X^n_s, (X^n_s \cdot\nabla)\psi_{z,j} \rangle
-\langle\tau(X^n_s), e (\psi_{z,j})\rangle,\nonumber
\\[-8pt]\\[-8pt]
      W^{z,j}_t &=&\langle W_t, \psi_{z,j} \rangle.\nonumber
\end{eqnarray}
Let $W_\cdot$ and $\xi$ be as above. We then define
\begin{eqnarray*}
\mathcal{G}^{\xi,W}_t & = &\sigma( \xi, W_s,   s \le t), \qquad    0
\le t <\infty,
 \qquad     \mathcal{G}^{\xi,W}_\infty = \sigma \biggl( \bigcup_{t \ge
0}\mathcal{G}^{\xi,W}_t\biggr),
\\
\mathcal{N}^{\xi,W} & = &\{ N \subset\Omega   ;   \exists
\widetilde{N} \in\mathcal{G}^{\xi
,W}_\infty,
N \subset\widetilde{N},   P^W(\widetilde{N})=0 \}
\end{eqnarray*}
and
\begin{equation}\label{cF^W_t}
\mathcal{F}^{\xi,W}_t = \sigma( \mathcal{G}^{\xi,W}_t \cup
\mathcal{N}^{\xi,W} ),
  \qquad  0 \le t <\infty.
\end{equation}
In what follows, expectation with respect to the measure $P^W$ will be
denoted by~$E^W [  \cdot  ]$.

\begin{theorem}\label{Thm.Gar}
Let $W_\cdot$, $\xi$ and $\mathcal{F}^{\xi,W}_t$ be as above.
Then for each $n =1,2,\ldots$
there exists a unique process $X^n_\cdot$ such that:
\begin{enumerate}[(a)]
\item[(a)] $X^n_t$ is $\mathcal{F}^{\xi,W}_t$-measurable for all $t
\ge0$;
\item[(b)] $(\ref{sde_Gar^n})$ is satisfied;
\item[(c)] For any $T>0$,
\begin{eqnarray}\label{ener=}
E^W \biggl[ \| X^n_T \|_2^2+2\int^T_0\langle e(X^n_t),\tau(X^n_t)
\rangle\, dt
\biggr]
& = & E^W [ \| X^n_0 \|_2^2 ] +\operatorname{tr}(\Gamma\mathcal{P}_n )T,
\\ \label{apriori_p}
E^W \biggl[ \| X^n_T \|_2^2+\frac{1 }{ C}\int^T_0\| X^n_t \|
_{p,1}^p\,dt\biggr]
& \le& m_0 +\bigl(C+\operatorname{tr}(\Gamma)\bigr)T <\infty,
\end{eqnarray}
where $C=C(d, p) \in(0,\infty)$.
\end{enumerate}
Suppose, in addition, that $p \ge\frac{2d }{ d+2}$,
where $p$ is from $(\ref{t(v)})$.
Then, for any $T>0$,
\begin{equation}\label{apriori}
E^W \biggl[ \sup_{t \le T}\| X^n_t \|_2^2 +\int^T_0\| X^n_t \|
_{p,1}^p\,dt\biggr] \le
(1+T)C' <\infty,
\end{equation}
where $C'=C'(d, p,\Gamma,m_0) \in(0,\infty)$.
\end{theorem}

\begin{pf}
We fix the accuracy $n$ of the approximation introduced above
and suppress the superscript ``$n$'' from the notation $X=X^n$.
We write the summation over $z \in[-n,n]^d$ and $j=1,\ldots,2d-2$ simply by
$\sum_{z, j}$.
Since $v \mapsto\mathcal{P}_n b (v)\dvtx \mathcal{V}_n \rightarrow
\mathcal{V}_n$ is locally
Lipschitz continuous [see (\ref{b^zj})] and
\renewcommand{\theequation}{1}
\setcounter{equation}{0}
\begin{equation}
\langle v, b(v) \rangle
\stackrel{\mbox{\scriptsize(\ref{uu=0})}}{=}
-\langle e (v), \tau(v) \rangle\le
C-\frac{1 }{ C}\| v \|_{p,1}^p,
\end{equation}
where we have used \cite{MNRR96}, formula (1.11), page~196, and formula $(1.20)_2,$ page~198,
to see the second inequality.
This implies that there exists a unique process~$X_\cdot$ with
the properties (a)--(b) above, as can be seen from
standard existence and uniqueness results for the SDE,
for example, \cite{IW89}, Theorem 2.4, page~177, and Theorem~3.1, pages~178--179;
cf. the remark after the proof.
Note that for $\alpha=0,1,2,\ldots,$
\[
\| \nabla^\alpha v\|_2^2 =\langle v, (-\Delta)^\alpha v \rangle=\sum
_{z,j} (-4\pi^2
|z|^2)^\alpha\langle v, \psi_{z,j} \rangle^2,
     \qquad v \in\mathcal{V}_n.
\]
On the other hand, we have by It\^o's formula that
\[
|X^{z,j}_t|^2=|X^{z,j}_0|^2+2\int^t_0X^{z,j}_s\,dW^{z,j}_s+2\int
^t_0X^{z,j}_sb^{z,j}_s(X_s)\,ds
+\langle\psi_{z,j},\Gamma\psi_{z,j}\rangle t.
\]
Therefore,
\renewcommand{\theequation}{\arabic{section}.\arabic{equation}}
\setcounter{equation}{10}
\begin{eqnarray}\label{Ito^2}
\| \nabla^\alpha X_t \|_2^2&=&\| \nabla^\alpha X_0 \|_2^2+ 2M_t
+2\int^t_0 \langle(-\Delta)^\alpha X_s, b(X_s) \rangle\, ds\nonumber
\\[-8pt]\\[-8pt]
&&{}+\operatorname{tr} (\Gamma
(-\Delta)^\alpha\mathcal{P}_n)t,\nonumber
\end{eqnarray}
where
\begin{equation}\label{Ito^2M}
M_t = \sum_{z,j}\int^t_0 (-\Delta)^\alpha X^{z,j}_s \,d W^{z,j}_s.
\end{equation}
Here we will use (\ref{Ito^2}) only for $\alpha=0$. The case $\alpha
=1$ will
be used in the proof of Lemma~\ref{Lem.DX} later on.
By (\ref{Ito^2}) with $\alpha=0$,
\renewcommand{\theequation}{2}
\setcounter{equation}{0}
\begin{equation}
\| X_t \|_2^2+\frac{2 }{ C}\int^t_0 \| X_s\|_{p,1}^p\, ds
\le\| X_0 \|_2^2+2M_t +\bigl(C +\operatorname{tr} (\Gamma) \bigr)t,
\end{equation}
where $M_t$ in (2) is defined by (\ref{Ito^2M}) with $\alpha=0$.
Since it is not difficult to see that the above $M_t$ is a martingale
(cf. \cite{Fl08}, proof of (10), page~60),
we get (\ref{ener=}) by taking expectation of the equality (\ref{Ito^2}).
Similarly, we obtain (\ref{apriori_p})
by taking expectation\vadjust{\goodbreak} of the inequality (2).
To see (\ref{apriori}), it is enough to show that
there exists $\delta\in(0,1]$ such that\vspace*{-2pt}
\renewcommand{\theequation}{3}
\setcounter{equation}{0}
\begin{equation}
 E^W \Bigl[ \sup_{t \le T}\| X_t \|_2^2 \Bigr]
\le(1+T)C+C E^W \biggl[ \biggl( \int^T_0\| X_t \|_{p,1}^p\,dt
\biggr)^\delta\biggr] .\vspace*{-2pt}
\end{equation}

\noindent
To see this, we start with a
bound on the quadratic variation of the martingale $M_\cdot$,\vspace*{-2pt}
\renewcommand{\theequation}{4}
\setcounter{equation}{0}
\begin{equation}
\langle M \rangle_t = \int^t_0 \langle\Gamma X_s,
X_s \rangle \,ds
\le\| \Gamma\|_{2 \rightarrow2} \int^t_0 \| X_s \|_2^2\,ds,\vspace*{-2pt}
\end{equation}

\noindent
where $\| \Gamma\|_{2 \rightarrow2}$ denotes the operator norm of
$\Gamma\dvtx V_{2,0} \rightarrow V_{2,0}$.
We now recall the Burkholder--Davis--Gundy inequality
(\cite{IW89}, Theorem 3.1, page 110),\vspace*{-2pt}
\renewcommand{\theequation}{5}
\setcounter{equation}{0}
\begin{equation}
 E^W \Bigl[ \sup_{t \le T}| M_t |^q \Bigr]
\le C E^W [ \langle M \rangle_T^{q/2}]    \quad   \mbox{for } q
\in(0,\infty).\vspace*{-2pt}
\end{equation}

\noindent
We then observe that\vspace*{-2pt}
\renewcommand{\theequation}{6}
\setcounter{equation}{0}
\begin{eqnarray}
 E^W \Bigl[ \sup_{t \le T}\| X_t \|_2^2 \Bigr]
& \stackrel{\mbox{\scriptsize(2)}}{\le} &
(1+T)C +2E^W \Bigl[ \sup_{t \le T}| M_t | \Bigr]\nonumber
\\[-8pt]\\[-8pt]
& \stackrel{\mbox{\scriptsize(4)--(5)}}{\le} &
(1+T)C +C'E^W \biggl[ \biggl( \int^T_0 \| X_s \|_2^2\,ds
\biggr)^{1/2}\biggr].\nonumber\vspace*{-2pt}
\end{eqnarray}

\noindent
This proves (3) for $p \ge2$. We assume $p<2$ in what follows.
We have\vspace*{-2pt}
\[
e_\ell\stackrel{\rm def.}{=}\inf\{ t   ;   \| X_t \|_2 \ge\ell\}
\nearrow\infty,    \qquad \mbox{as $\ell\nearrow\infty$,}\vspace*{-2pt}
\]
since the process $X_t$ does not explode.
On the other hand, it is clear that
the following variant of (6) is true:\vspace*{-2pt}
\renewcommand{\theequation}{6$'$}
\setcounter{equation}{0}
\begin{eqnarray}
E^W \Bigl[ \sup_{t \le T \wedge e_\ell}\| X_t \|_2^2 \Bigr]
\le
(1+T)C +CE^W \biggl[ \biggl( \int^{T \wedge e_\ell}_0
\| X_s \|_2^2\,ds \biggr)^{1/2}\biggr].\vspace*{-2pt}
\end{eqnarray}

\noindent
We have by Sobolev embedding that for $v \in V_{p,1}$,\vspace*{-2pt}
\renewcommand{\theequation}{7}
\setcounter{equation}{0}
\begin{eqnarray}
\| v \|_2 \le C \|v \|_{p,1},  \qquad   \mbox{since } p \ge\frac{2d }{ d+2}.\vspace*{-2pt}
\end{eqnarray}

\noindent
Let $\varepsilon>0$, $r=\frac{4 }{2-p} \in(4,\infty)$ and
$r'=\frac{r }{ r-1} =\frac{4 }{2+p} \in(1,4/3)$. Then,\vspace*{-2pt}
\renewcommand{\theequation}{8}
\setcounter{equation}{0}
\begin{eqnarray}
&&\biggl( \int^{T \wedge e_\ell}_0 \| X_s \|_2^2\,ds
\biggr)^{1/2}\nonumber
\\[-1pt]
&&\hspace*{6,5pt}\qquad \le \sup_{s \le T \wedge e_\ell}\|
X_s\|_2^{{(2-p)/2}}
\biggl( \int^{T \wedge e_\ell}_0 \| X_s\|_2^p \,ds \biggr)^{1/2}\nonumber
\\[-8pt]\\[-8pt]
&&\hspace*{6,5pt}\qquad \stackrel{\mbox{\scriptsize(7)}}{\le}
 C\sup_{s \le T \wedge e_\ell}\| X_s\|_2^{{(2-p )/2}}
\biggl( \int^{T \wedge e_\ell}_0 \| X_s\|_{p,1}^p\, ds \biggr)^{1/2}\nonumber
\\[-1pt]
&&\qquad \stackrel{\mbox{\scriptsize Young}}{\le}
\frac{\varepsilon^r C }{ r}\sup_{s \le T \wedge e_\ell
}\| X_s\|_2^2
+\frac{\varepsilon^{-r'} C }{ r'}
\biggl( \int^{T \wedge e_\ell}_0 \| X_s\|_{p,1}^p \,ds \biggr)^{{2
/(2+p)}}.\nonumber
\end{eqnarray}
Since $E^W [ \sup_{t \le T \wedge e_\ell}\| X_t \|_2^2 ]
\le
\ell^2 <\infty$,
we have by (6) and (8) that
\[
E^W \Bigl[ \sup_{t \le T \wedge e_\ell}\| X_t \|_2^2 \Bigr]
\le(1+T)C+C E^W \biggl[ \biggl( \int^{T \wedge e_\ell}_0
\| X_t \|_{p,1}^p\,dt \biggr)^{{2 /(2+p)}} \biggr].
\]
Letting $\ell\nearrow\infty$, we obtain (3).
\end{pf}

\begin{remark*}
Unfortunately, the SDE (\ref{sde_Gar^n}) does not
satisfy the condition (2.18) imposed in the existence theorem
(\cite{IW89}, Theorem 2.4, page~177). However,
we easily see from the proof of the existence theorem
that (2.18) there can be replaced by
\[
\| \sigma(x)\|^2+x \cdot b(x) \le K(1+|x|^2).
\]
We have applied \cite{IW89}, Theorem 2.4, page~177, with this
modification.
\end{remark*}

\subsection{Further a priori bounds}

We first prove the following general estimates which apply both
to the weak solution $X$ to (\ref{SNS_def}) and to the
unique solution to~(\ref{sde_Gar^n}).
%
\begin{lemma}\label{Lem.intb(X)}
Let $T>0$ and $X=(X_t)_{t \ge0}$ be a process on a probability space
$(\Omega, \mathcal{F}, P)$ such that
\[
X \in L_p ([0,T] \rightarrow V_{p,1})
\cap L_{\infty} ([0,T] \rightarrow V_{2,0}),   \qquad  \mbox{a.s.}
\]
and
\[
A_T=E \biggl[ \int^T_0 \| X_s \|_{p,1}^{p}\,ds \biggr]<\infty,
\qquad
B_T=E \Bigl[ \sup_{s \in[0,T]}\| X_s \|_2^2 \Bigr]<\infty.
\]
\begin{longlist}
\item[\textup{(a)}]
For $p \in[\frac{2d }{ d+2}, \infty)$,
\renewcommand{\theequation}{\arabic{section}.\arabic{equation}}
\setcounter{equation}{12}
\begin{equation}\label{intb(X)a}
E \biggl[
\biggl( \int^T_0 \| (X_s \cdot\nabla)X_s \|_{p',-\beta(p,1)}^p\,ds
\biggr)^\delta\biggr]
\le CA_T^\delta B_T^{1-\delta} <\infty,
\end{equation}
where $\delta=\frac{p }{ p+2}$, $p'=\frac{p }{ p-1}$,
$\beta(p,1)$ is defined by $(\ref{a_4=1})$
and $C=C(d, p) \in(0,\infty)$.
\item[\textup{(b)}]
\begin{equation}\label{intb(X)b}
E \biggl[ \int^T_0 \| \operatorname{div}\tau(X_s) \|_{p',-1}^{p'}
\,ds\biggr]
\le(T+A_T)C'<\infty,
\end{equation}
where $C'=C'(p,\nu) \in(0,\infty)$.
\end{longlist}
\end{lemma}

\begin{pf} (a)
We have by (\ref{vwinV}) that
\renewcommand{\theequation}{1}
\setcounter{equation}{0}
\begin{equation}
\| (v \cdot\nabla) v \|_{p', -\beta(p,1)}
\le C \| v\|_{p,1}\| v\|_{2}
\qquad\mbox{for } v \in V_{p,1} \cap V_{2,0}.
\end{equation}
We then use (1) to see that
%
\begin{eqnarray*}
I &\stackrel{\rm def.}{=}&
\int^T_0 \| (X_s \cdot\nabla)X_s \|_{p',-\beta(p,1)}^{p}\,ds
\stackrel{\mbox{\scriptsize(1)}}{\le} C
\int^T_0 \| X_s \|_{p,1}^{p}\| X_s \|_2^{ p}\,ds
 \\
& \le& C \sup_{s \in[0,T]}\| X_s \|_2^{p}
\int^T_0 \| X_s \|_{p,1}^{p}\,ds.
\end{eqnarray*}
%
Finally, noting that $\frac{p\delta}{1-\delta}=2$, we
conclude that
\begin{eqnarray*}
E[I^\delta]
& \le&
C E \biggl[ \sup_{s \in[0,T]}\| X_s \|_2^{p \delta}
\biggl( \int^T_0 \| X_s \|_{p,1}^{p}\,ds\biggr)^\delta\biggr]
\\
& \le&
C E \Bigl[ \sup_{s \in[0,T]}\| X_s \|_2^2 \Bigr]^{1-\delta}
E \biggl[ \int^T_0 \| X_s \|_{p,1}^{p}\,ds \biggr]^\delta
= CB_T^{1-\delta}A_T^\delta.
\end{eqnarray*}

(b)
\[
\| \operatorname{div}\tau(X_s) \|_{p',-1}
\stackrel{\mbox{\scriptsize(\ref{b_2<})}}{\le} C
\bigl(1+\|e (X_s)\|_p\bigr)^{p-1}
\]
which implies that
\[
\| \operatorname{div}\tau(X_s) \|_{p',-1}^{p'}
\le C+C\|e (X_s)\|_p^p
\]
and hence, that
\begin{eqnarray*}
&&E \biggl[\int^T_0 \| \operatorname{div}\tau(X_s) \|
_{p',-1}^{p'}\,ds\biggr]
 \\
&&\qquad \le
CT +C E \biggl[ \int^T_0 \| e (X_s)\|_p^p\,ds \biggr]\le(T+A_T)C.
\end{eqnarray*}
\upqed\end{pf}

Let $X^n=(X^n_t)_{t \ge0} \in\mathcal{V}$ be the unique
solution of (\ref{sde_Gar^n})
for the Galerkin approximation.

\begin{lemma}\label{Lem.3.34}
Suppose $(\ref{<p<})$.
Then, there exist $\tilde{p} \in(1,p)$ and $\tilde{\alpha}
\in(1,\infty)$
such that for each $T>0$
\renewcommand{\theequation}{\arabic{section}.\arabic{equation}}
\setcounter{equation}{14}
\begin{equation}\label{3.34}
E^W \biggl[ \int^T_0\|X^n_t\|_{\tilde{p}, \tilde{\alpha
}}^{\tilde{p}}\,dt
\biggr] \le C_T <\infty,
\end{equation}
where the constant $C_T$ is independent of $n$.
\end{lemma}

We will have slightly better than
is stated in Lemma~\ref{Lem.3.34} in the course of the proof.
For (i) $d=2$ and $p \ge2$ and (ii) $d \ge3$ and $p >p_3 (d)$,
we have that
\begin{equation}\label{3.34''}
E^W \biggl[ \int^T_0\|\Delta X^n_t\|_2^{{2p/( p+2\lambda)}}\,dt
\biggr] \le C_T <\infty,
\end{equation}
where $\lambda\ge0$ is defined by (\ref{lm}) below.
For $p < \frac{2d }{ d-2}$, we have that
\begin{equation}\label{3.34'}
E^W \biggl[ \int^T_0\|X^n_t\|_{p, \tilde{\alpha}}^{\tilde{p}}\,dt
\biggr] \le C_T <\infty
\end{equation}
for \textit{any} $\tilde{p} \in(1,p)$ with some $\tilde{\alpha
}=\tilde{\alpha} (\tilde
{p}) >1$.

The rest of this section is devoted to the proof of Lemma~\ref{Lem.3.34}.
We suppress the superscript $n$ from the notation.
We write the summation over $z \in[-n,n]^d$ and $j=1,\ldots ,2d-2$ simply by
$\sum_{z, j}$.
We first
establish the following bounds.

\begin{lemma}\label{Lem.DX}
Suppose that $p \in(\frac{3d-4 }{ d},\infty)$ if $d \ge3$ and let
\begin{eqnarray}\label{lm}
\lambda&= & \cases{
0, &\quad \mbox{if $d = 2$},
\cr
\dfrac{2(3-p)^+}{ dp -3d+4}, &\quad \mbox{if $d \ge3$},
}
\\
\eqntext{\mbox{cf. \cite{MNRR96}, formula $(3.47)$,
page~236,}}
\\\label{cJ_t}
\mathcal{J}_t &= & \cases{
\dfrac{\| \Delta X_t \|_2^2 }{(1+\| \nabla X_t \|
^2_2)^\lambda},
&\quad \mbox{if $p \ge2$},
 \cr
\dfrac{\| \Delta X_t \|_p^2 }{
(1+\| \nabla X_t \|^2_2)^\lambda(1+\| \nabla X_t \|_p)^{2-p} },
&\quad \mbox{if $1<p < 2$}.
}
\end{eqnarray}
Then, for any $T>0$,
\begin{equation}\label{DX}
E^W \biggl[
\int^T_0\mathcal{J}_t\,dt \biggr] \le C_T <\infty,
\end{equation}
where $C_T=C(T, d, p, \Gamma, m_1 )$.
\end{lemma}

\begin{pf}
By (\ref{Ito^2}) with $\alpha=1$,
\renewcommand{\theequation}{1}
\setcounter{equation}{0}
\begin{equation}
\frac{1}{2}\| \nabla X_t \|_2^2 =\frac{1}{2}\| \nabla X_0 \|_2^2 +M_t
+\int^t_0 K_s\,ds,
\end{equation}
where
\[
M_t
= -\sum_{z,j}\int^t_0\Delta X^{z,j}_s \,d W^{z,j}_s,\qquad
K_s = \langle-\Delta X_s, b (X_s) \rangle+\frac{1}{2}\operatorname{tr}
(-\Gamma\Delta\mathcal{P}_n ).
\]

\textit{Step} $1$. We will prove that
\renewcommand{\theequation}{2}
\setcounter{equation}{0}
\begin{equation}
K_s +c_1 \mathcal{I}_s \le
\cases{ 0, &\quad \mbox{if $d =2$}, \cr
C_1(1+\| \nabla X_t \|^2_2)^\lambda(1+\| \nabla X_t \|_p)^p, &\quad \mbox
{if $d
\ge3$},
}
\end{equation}
where $c_1,C_1 \in(0,\infty)$ are constants and
\[
\mathcal{I}_s =\int_{{\mathbb{T}}^d}\bigl(1+|e (X_s)|^2\bigr)^{{(p-2 )/
2}}|\nabla e (X_s)|^2.
\]
%
To show (2), note that
\[
\langle-\Delta X_s, b (X_s) \rangle
= \langle-\Delta X_s, (X_s \cdot\nabla) X_s \rangle
-\langle\tau(X_s), e(-\Delta X_s) \rangle.
\]
We see from the
argument in \cite{MNRR96}, proof of (3.19), page~225, that
\renewcommand{\theequation}{3}
\setcounter{equation}{0}
\begin{equation}
\langle\tau(X_s), e(-\Delta X_s) \rangle\ge2c_1 \mathcal{I}_s.
\end{equation}
On the other hand, we have by integration by parts
and H\"older's inequality that
\[
\langle-\Delta X_s, (X_s \cdot\nabla) X_s \rangle
= \sum_{i,j,k}\int_{{\mathbb{T}}^d}\partial_k X^j_s \,\partial_j
X^i_s \,\partial_k X^i_s
\le\| \nabla X_s \|_3^3,
\]
where $X^j_s=\sum_{z \in[-n,n]^d}X^{z,j}_s \psi_{z,j}$.
It is also well known that the inner product on the LHS vanishes if $d=2$
(\cite{MNRR96}, formula (3.20), page~225).
By the argument in \cite{MNRR96}, proof of (3.46), pages~234--235
(this is where the choice of $\lambda$ is used),
we get
\[
\| \nabla X_s \|_3^3 \le C_1(1+\| \nabla X_t \|^2_2)^\lambda(1+\|
\nabla X_t \|
_p)^p+c_1\mathcal{I}_s.
\]
These imply that
\renewcommand{\theequation}{4}
\setcounter{equation}{0}
\begin{eqnarray}
&&\langle-\Delta X_s, (X_s \cdot\nabla) X_s \rangle\nonumber
\\[-8pt]\\[-8pt]
&&\qquad{}\times\cases{ =0, & \quad\mbox{if $d =2$}, \cr
\le C_1(1+\| \nabla X_t \|^2_2)^\lambda(1+\| \nabla X_t \|
_p)^p+c_1\mathcal{I}_s, &
\quad\mbox{if $d \ge3$}.
}\nonumber
\end{eqnarray}
We get (2) by (3)--(4).

\textit{Step} $2$.
Proof of (\ref{DX}).
By \cite{MNRR96}, formulas (3.25) and (3.26), page~227,
$\mathcal{J}_t$ and $\mathcal{I}_t$ are related as
\[
\mathcal{J}_t \le C\frac{\mathcal{I}_t }{(1+\| \nabla X_t \|
_2^2)^\lambda}.
\]
Therefore, it is enough to prove that
\renewcommand{\theequation}{5}
\setcounter{equation}{0}
\begin{eqnarray}
E^W \biggl[ \int^t_0\frac{\mathcal{I}_s\, ds}{(1+\| \nabla X_s \|_2^2
)^\lambda}\biggr]
\le C_T <\infty ,
\end{eqnarray}
where $C_T =C(T, d, p, \Gamma, m_0, m_1) \in(0,\infty)$.

To see this, we introduce the following concave function of $x \ge0$:
\[
f(x)=\cases{
\dfrac{1 }{1-\lambda}(1+x)^{1-\lambda}, &\quad \mbox{if $\lambda\neq1$}, \cr
\ln(1+x), &\quad \mbox{if $\lambda= 1$.}
}
\]
Then we have by (1) and It\^o's formula that
\[
f(\| \nabla X_t \|_2^2 )
\le f(\| \nabla X_0 \|_2^2 )
+\int^t_0\frac{dM_s }{(1+\| \nabla X_s \|_2^2 )^\lambda}
+2\int^t_0 \frac{K_s\, ds }{(1+\| \nabla X_s \|_2^2 )^\lambda},
\]
where we have omitted the term with $f'' \le0$.
Moreover, by (2)
\begin{eqnarray*}
\frac{K_s }{(1+\| \nabla X_s \|_2^2 )^\lambda} &\le&
-\frac{c_1 \mathcal{I}_s }{(1+\| \nabla X_s \|_2^2 )^\lambda}
+C_1 (1+\| \nabla X_s \|_p )^p,
\\
0 &\le& f(x) \le C_2 (1+x) \qquad  \mbox{if $\lambda\in[0,1]$}
\end{eqnarray*}
{and}
\[
-\frac{1 }{\lambda-1} \le f(x) \le0  \qquad \mbox{if $\lambda>1$}.
\]
Putting these together, we get
\begin{eqnarray*}
&&-C_3+2c_1E^W \biggl[
\int^t_0\frac{ \mathcal{I}_s \,ds }{(1+\| \nabla X_s \|_2^2 )^\lambda
}\biggr]
\\
&&\hspace*{5pt}\qquad \le C_2 (1+E[\| \nabla X_0\|_2^2])
+ C_1 E^W \biggl[ \int^t_0 (1+\| \nabla X_s \|_p )^p \,ds \biggr]
\\
&&\qquad \stackrel{\mbox{\scriptsize(\ref{apriori})}}{\le}
C(T, d, p, \Gamma, m_0, m_1) <\infty,
\end{eqnarray*}
where $C_3=0$ if $\lambda\in(0,1]$
and $C_3=\frac{1 }{\lambda-1}$ if $\lambda>1$. This proves (5).
\end{pf}

\begin{pf*}{Proof of Lemma~\ref{Lem.3.34}}
We note that
\begin{eqnarray*}
 p_1 (d) &<& p_3 (d) < p_2 (d)  \qquad   \mbox{for $d \le8$},
 \\
 p_1 (9)&=&2.555\ldots<p_2 (9)=2.5714\ldots<p_3(9)=2.620\ldots,
 \\
 p_2 (d) &<& p_1 (d)\qquad    \mbox{for $d \ge10$}.
\end{eqnarray*}
Thus, condition (\ref{<p<}) takes the following form in any $d \ge2$:
\renewcommand{\theequation}{\arabic{section}.\arabic{equation}}
\setcounter{equation}{20}
\begin{equation}\label{<p<2}
p \in( p_1 (d),p_2 (d)) \cup
(p_3 (d), \infty ).
\end{equation}
We consider the following four cases separately:
\begin{longlist}
\item[\textit{Case} $1$.] $d =2$ and $p \ge2$;

\item[\textit{Case} $2$.] $d \ge3$ and $p >p_3 (d)$;

\item[\textit{Case} $3$.] $p \in(p_1(d),p_2 (d))$ and $p \ge2$;

\item[\textit{Case} $4$.] $p \in(p_1 (d), 2)$ (this case appears only
if $d =2,3$).
\end{longlist}

The first two cases cover the interval $(p_3 (d), \infty)$ in (\ref
{<p<2}). [Note that $p_3 (2)=2$,
while the last two cases cover the interval $(p_1 (d), p_2 (d)$.]

\textit{Case} $1$.
By (\ref{DX}),
(\ref{3.34}) has already been shown with $\tilde{p}=\tilde
{\alpha} =2$.\vadjust{\goodbreak}

\textit{Case} $2$. Note that $p >p_3 (d)>2$ and that $\beta
\stackrel
{\rm def.}{=}\frac{p }{ p+2\lambda} >1/2$.
We prove (\ref{3.34''}).
Since $\lambda\beta=\frac{p }{2}(1-\beta)$,
\renewcommand{\theequation}{1}
\setcounter{equation}{0}
\begin{eqnarray}
 &&E^W \biggl[ \int^T_0 \| \Delta X_s \|_2 ^{2 \beta}ds\biggr]\nonumber
\\
&&\hspace*{19pt}\qquad=  E^W \biggl[ \int^T_0 \mathcal{J}_s^\beta
(1+\| \nabla X_s \|_2^2 )^{ \lambda\beta}\,ds \biggr]\nonumber
\\[-8pt]\\[-8pt]
&&\qquad\hspace*{1pt} \stackrel{\beta+(1 -\beta)=1}{\le}
 E^W \biggl[ \int^T_0\mathcal{J}_s \,ds\biggr]^\beta
E^W \biggl[ \int^T_0(1+\| \nabla X_s \|_2^2)^{{p /2}}
\,ds \biggr]^{1 -\beta}\nonumber
 \\
&&\ \qquad \stackrel{\mbox{\scriptsize(\ref{apriori}), (\ref{DX})}}{\le}
C_T <\infty,\nonumber
\end{eqnarray}
where we used (\ref{DX}) for $p \ge2$.

\textit{Case} $3$.
We prove (\ref{3.34'}) for given $\tilde{p} \in(1,p)$
with some $\tilde{\alpha}=\tilde{\alpha} (\tilde{p})
\in(1,2)$.
Let $\beta=\frac{p }{ p+2\lambda} \in(0,1)$.
Then the bound (1) from case 2 is still valid, although it
may no longer be the case that $2\beta>1$ here.
On the other hand, it is not difficult to see via
the interpolation and the Sobolev imbedding that
for any $\tilde{p} \in(1,p)$, there exist $\tilde{\alpha}
\in(1,2)$ and
$\theta\in(0,1)$
such that
\[
\int^T_0 \| X_s \|_{p,\tilde{\alpha} }^{\tilde{p}}\,ds
\le C \biggl( \int^T_0 \| X_s \|_{p,1}^p\,ds \biggr)^\theta
\biggl( \int^T_0 \| X_s \|_{2,2}^{2\beta}\,ds \biggr)^{1-\theta};
\]
cf. \cite{MNRR96}, proof of (3.58), page~238. This is where the restriction
$p <\frac{2d }{ d-2}$ is necessary.
Thus,
\begin{eqnarray*}
E^W \biggl[ \int^T_0 \| X_s \|_{p,\tilde{\alpha} }^{\tilde
{p}}\,ds \biggr]
&\le& C E^W \biggl[ \int^T_0 \| X_s \|_{p,1}^p\,ds \biggr]^\theta
E^W \biggl[ \int^T_0 \| X_s \|_{2,2}^{2\beta}\,ds \biggr]^{1-\theta}
\\
&\stackrel{\mbox{\scriptsize(\ref{apriori}), (1)}}{\le}& C_T <\infty.
\end{eqnarray*}
%

\textit{Case} $4$.
We prove (\ref{3.34'}) for given $\tilde{p} \in(1,p)$
and with some $\tilde{\alpha}=\tilde{\alpha} (\tilde
{p}) \in(1,2)$.
We recall that $p >\frac{3d }{ d+2}$ and set
\[
\beta=\frac{((d+2)p-3d)p }{2((d+5)p-3d-p^2)} \in\biggl(0,\frac{1}{2}\biggr).
\]
Then,
\renewcommand{\theequation}{2}
\setcounter{equation}{0}
\begin{eqnarray}
\rho\stackrel{\rm def.}{=} \frac{(2-p)d \lambda}{2(1-\beta)p }\in[0,1)
   \quad \mbox{and} \quad   \frac{(2-p)\beta}{1-\beta} \in(0,p).
\end{eqnarray}
As a result of applications of H\"older's inequality,
the interpolation and the Sobolev imbedding
(cf. \cite{MNRR96}, formulas  (3.60)--(3.63), pages 239--240), we arrive at the
following bound:
\renewcommand{\theequation}{3}
\setcounter{equation}{0}
\begin{eqnarray}
\int^T_0 \| \Delta X_s \|_p^{2\beta}\,ds
\le C \biggl( \int^T_0 \mathcal{J}_s \,ds \biggr)^\beta
( I_1+I_2)^{1-\beta},
\end{eqnarray}
where
\begin{eqnarray*}
I_1 &=&\int^T_0 (1+\| \nabla X_s \|_p)^{{(2-p)\beta/(1-\beta)}}\,ds,
   \\[-2pt]
    I_2&=&
\biggl( \int^T_0 \| \Delta X_s \|_p^{2\beta}\,ds \biggr)^\rho
\biggl( \int^T_0 \| \nabla X_s \|_p^p\,ds \biggr)^{1-\rho}.
\end{eqnarray*}
We first prove that
\renewcommand{\theequation}{4}
\setcounter{equation}{0}
\begin{eqnarray}
E^W \biggl[ \int^T_0 \| \Delta X_s \|_p^{2\beta}\,ds \biggr]
\le C_T <\infty.
\end{eqnarray}
We first assume $d=3$, where $\rho>0$. Let $r= \frac{1 }{\rho}\in
(1,\infty)$
and $r'=\frac{r }{ r-1}=\frac{1 }{1-\rho} \in(1,\infty)$.
Then, for $\varepsilon>0$,
\begin{eqnarray*}
E^W \biggl[ \int^T_0 \| \Delta X_s \|_p^{2\beta}\,ds \biggr]
& \stackrel{\mbox{\scriptsize(3)}}{\le}&
C E^W \biggl[ \biggl( \int^T_0 \mathcal{J}_s \,ds \biggr)^{\beta}
( I_1+I_2)^{1-\beta} \biggr] \nonumber
\\[-2pt]
& \stackrel{\beta+(1 -\beta)=1}{\le} &
CE^W \biggl[ \int^T_0 \mathcal{J}_s\, ds \biggr]^\beta E^W [
I_1+I_2]^{1-\beta
} \nonumber
\\[-2pt]
& \stackrel{\mbox{\scriptsize(\ref{DX})}}{\le}&
C_TE[ 1+ I_1+I_2],\nonumber
\\[-2pt]
E^W [ I_1]
& \stackrel{\mbox{\scriptsize(\ref{apriori}),(2)}}{\le}&
C_T <\infty, \nonumber
\\[-2pt]
E^W [ I_2 ]
& \stackrel{\mbox{\scriptsize Young}}{\le} &
\frac{\varepsilon^r }{ r}
E^W \biggl[
\int^T_0 \| \Delta X_s \|_p^{2\beta}\,ds \biggr]
\\[-2pt]
&&{}+\frac{\varepsilon^{-r'}
}{ r'}
E^W \biggl[
\int^T_0 \| \nabla X_s \|_p^p\,ds \biggr] \nonumber
\\[-2pt]
& \stackrel{\mbox{\scriptsize(\ref{apriori}) }}{\le} &
\frac{\varepsilon^r }{ r}
E^W \biggl[
\int^T_0 \| \Delta X_s \|_p^{2\beta}\,ds \biggr] +C_T. \nonumber
\end{eqnarray*}
Putting things together, with $\varepsilon$ small enough, we arrive at
(4) for $d=3$.
If $d=2$ and hence, $\rho=0$, then we have $E^W [ I_2 ]
\le
C_T$ directly
from (\ref{apriori}). Therefore, the proof of (4) is even easier than
the above.

We finally turn to (\ref{3.34}).
It is not difficult to see via
the interpolation
(cf. \cite{MNRR96}, proof of (3.65), pages~240--241) that
for any $\tilde{p} \in(1,p)$, there exist $\tilde{\alpha}
\in(1,2)$ and
$\theta\in(0,1)$
such that
\[
\int^T_0 \| X_s \|_{p,\tilde{\alpha} }^{\tilde{p}}\,ds
\le C \biggl( \int^T_0 \| X_s \|_{p,1}^p\,ds \biggr)^\theta
\biggl( \int^T_0 \| X_s \|_{p,2}^{2\beta}\,ds \biggr)^{1-\theta}.
\]
Thus,
\begin{eqnarray*}
E^W \biggl[ \int^T_0 \| X_s \|_{p,\tilde{\alpha} }^{\tilde
{p}}\,ds \biggr]
&\le& C E^W \biggl[ \int^T_0 \| X_s \|_{p,1}^p\,ds \biggr]^\theta
E^W \biggl[ \int^T_0 \| X_s \|_{p,2}^{2\beta}\,ds \biggr]^{1-\theta}
\\
&\stackrel{\mbox{\scriptsize(\ref{apriori}),(4)}}{\le}& C_T
<\infty.
\end{eqnarray*}
\upqed\end{pf*}

\subsection{Compact imbedding lemmas}

We will need some compact imbedding
lemmas from \cite{FG95}.
We first introduce the following definition.

\begin{definition}\label{Def.Sob[0,T]}
Let $p \in[1,\infty)$, $T \in(0,\infty)$
and $E$ be a Banach space.
\begin{enumerate}[(a)]
\item[(a)]
We let
$L_{p,1}([0,T] \rightarrow E)$ denote the Sobolev space of
all $u \in L_p([0,T] \rightarrow E)$ such that
\[
u (t)=u (0)+\int^t_0 u'(s)\,ds\qquad     \mbox{for almost all $t \in[0,T]$}
\]
with some $ u(0) \in E$ and $u' (\cdot) \in L_p([0,T] \rightarrow E)$.
We endow the space $L_{p,1}([0,T] \rightarrow E)$ with the norm
$\| u \|_{L_{p,1}([0,T] \rightarrow E)}$ defined by
\[
\| u \|_{L_{p,1}([0,T] \rightarrow E)}^p
=\int^T_0\bigl(|u (t)|_E^p+|u'(t)|_E^p\bigr)\,dt.
\]
\item[(b)]
For $\alpha\in(0,1)$, we let
$L_{p,\alpha}([0,T] \rightarrow E)$ denote the Sobolev space of
all $u \in L_p([0,T] \rightarrow E) $ such that
\[
\int_{0 < s <t <T}\frac{|u(t)-u(s)|_E^p }{|t-s|^{1+\alpha p}}\,ds\,dt
<\infty.
\]
We endow the space $L_{p,\alpha}([0,T] \rightarrow E)$ with the norm
$\| u \|_{L_{p,\alpha}([0,T] \rightarrow E)}$ defined by
\[
\| u \|_{L_{p,\alpha}([0,T] \rightarrow E)}^p
=\int^T_0|u (t)|^p\,dt+
\int_{0 < s <t <T}\frac{|u(t)-u(s)|_E^p }{|t-s|^{1+\alpha p}}\,ds\,dt.
\]
\end{enumerate}
\end{definition}

To introduce the compact imbedding lemmas, we agree on the following standard
convention. Let $X$ be a vector space and $X_i \subset X$
be a subspace with the norm $\| \cdot\|_i$
($i=1,2$). Then we equip $X_0 \cap X_1$ and $X_0+X_1$, respectively,
with the norms
\begin{eqnarray*}
\| u\|_{X_0 \cap X_1}&=&\|u \|_0+\| u \|_1,
 \\
\| u\|_{X_0 + X_1}&=&\inf\{ \|u_0 \|_0+\| u_1 \|_1
  ;   u=u_0+u_1,   u_i \in X_i\}.
\end{eqnarray*}
The following lemmas will be used in Section \ref{conv_appro}.
%
\begin{lemma}[(\cite{FG95}, Theorem 2.2, page~370)]\label{Lem.cpt1}
Let:

\begin{itemize}
\item[$\bullet$]
$E_1,\ldots,E_n$ and $E$ be Banach spaces such that
each $E_i \stackrel{\rm compact}{\hookrightarrow} E$, $i=1,\ldots,n$.

\item[$\bullet$]
$p_1,\ldots,p_n \in(1,\infty)$, $\alpha_1,\ldots,\alpha_n >0$ are
such that $p_i\alpha_i >1$, $i=1,\ldots,n$.
\end{itemize}
Then, for any $T>0$,
\[
L_{p_1, \alpha_1} ([0,T] \rightarrow E_1)+\cdots+L_{p_n, \alpha_n}
([0,T] \rightarrow E_n)
\stackrel{\rm compact}{\hookrightarrow} C ([0,T] \rightarrow E).
\]
%
\end{lemma}

\begin{lemma}[(\cite{FG95}, Theorem 2.1, page~372)]\label{Lem.cpt2}
\it Let
\[
E_0 \stackrel{\rm compact}{\hookrightarrow} E \hookrightarrow E_1
\]
be Banach spaces such that the first embedding is compact
and $E_0,E_1$ are reflexive. Then, for any $p \in(1,\infty)$,
$\alpha\in(0,1)$ and $T>0$,
\[
L_p ([0,T] \rightarrow E_0) \cap L_{p,\alpha} ([0,T] \rightarrow E_1)
\stackrel{\rm compact}{\hookrightarrow} L_p ([0,T] \rightarrow E).
\]
\end{lemma}

\subsection{Convergence of the approximations} \label{conv_appro}

Let $X^n=(X^n_t)_{t \ge0} \in\mathcal{V}$ be the unique
solution to (\ref{sde_Gar^n})
for the Galerkin approximation. We write
\renewcommand{\theequation}{\arabic{section}.\arabic{equation}}
\setcounter{equation}{21}
\begin{equation}\label{p''}
p'=\frac{p }{ p-1}, \qquad    p''=p \wedge p'.
\end{equation}
Let $\beta(p,1)$ be defined by (\ref{a_4=1}) and
let $\tilde{p}>1$ be the one from Lemma~\ref{Lem.3.34}.
We may assume that $\tilde{p} \in(1, p'']$.
We also agree on the following standard convention. Let
$S$ be a set and $\rho_i $ be a metric on $S_i \subset S$ ($i=1,2$).
Then we tacitly consider the
metric $\rho_1+\rho_2$ on the set $S_1 \cap S_2$; cf. (\ref{space_tight}).

\begin{proposition}\label{Prop.tight}
Let $\beta> \beta(p,1)$. Then there exist
a process $X$ and a sequence $(\widetilde{X}^k)_{k \ge1}$ of processes
defined on a probability space $(\Omega, \mathcal{F}, P)$ such that
the following properties are satisfied:
\begin{enumerate}[(a)]
\item[(a)]
The process $X$ takes values in
\begin{equation}\label{space_tight}
C\bigl([0,\infty) \rightarrow V_{2 \wedge p',-\beta} \bigr)
\cap L_{\tilde{p},{\rm loc}}\bigl([0,\infty) \rightarrow V_{\tilde
{p},1} \bigr).
\end{equation}
\item[(b)]
For some sequence $n(k) \nearrow\infty$,
$\widetilde{X}^k$ has the same law as $X^{n(k)}$ and
\begin{equation}\label{tlX->X}
\lim_{k \rightarrow\infty}\widetilde{X}^k=X
\mbox{ in the metric space $(\ref{space_tight})$, $P$-a.s. }
\end{equation}
\end{enumerate}
\end{proposition}

\begin{remarks*} {(1)}
Due to Skorohod's representation theorem used in Lem\-ma~\ref
{Lem.tight0} below,
the probability space $(\Omega, \mathcal{F}, P)$ in the above
proposition may not
be the same as $(\Omega^W, \mathcal{F}^W, P^W)$,
where we have solved the SDE (\ref{sde_Gar^n}).

{(2)}
See (\ref{convp_1}) below for additional information on the
convergence (\ref{tlX->X}).
\end{remarks*}

We divide the Proposition~\ref{Prop.tight} into Lemmas~\ref
{Lem.tight1}--\ref{Lem.tight0}. To prepare the proofs of these lemmas,
we write (\ref{sde_Gar^n}) as
\begin{equation}\label{J+W}
X^n_t=X^n_0+I^n_t+J^n_t+W^n_t,
\end{equation}
with
\begin{eqnarray*}
I^n_t&=&\int^t_0 \mathcal{P}_n \bigl((X^n_s \cdot\nabla) X^n_s\bigr) \,ds,
\qquad
J^n_t=\int^t_0 \mathcal{P}_n (\operatorname{div}\tau(X^n_s))\, ds,
 \\
      W^n_t&=&\mathcal{P}_n W_t.
\end{eqnarray*}
It is elementary to obtain
the following regularity bound of the noise term $W^n_t$
\cite{Fl08}, Corollary 4.2, page~92:
for any $p \in[1,\infty)$, $\alpha\in[0,1/2)$ and $T>0$, there
exists $C_T=C_{\alpha, p, T} \in(0,\infty)$ such that
\begin{equation}\label{HBM1/2}
\sup_{n \ge0}E^W\bigl[\| W^n_\cdot\|_{L_{p,\alpha} ([0,T] \rightarrow
V_{2,0})}^p \bigr]
\le C_T{\rm tr}(\Gamma)^{p/2}.
\end{equation}
We will control $I^n_\cdot$ and $J^n_\cdot$ by
(\ref{intb(X)a}) and (\ref{intb(X)b}). However,
to be able to do so, we have to get rid of the projection
$\mathcal{P}_n$. This is the content of the following:

\begin{lemma}\label{Lem.I^nJ^n}
Let $T \in(0,\infty)$.
Then,
\begin{equation}\label{I^n}
\sup_{n \ge1}
E^W \bigl[ \| I^n_\cdot\|_{L_{p,1}([0,T] \rightarrow V_{p',-\beta
(p,1)})}^\gamma
\bigr]
\le C_T <\infty,
\end{equation}
where $\gamma=\frac{p^2 }{ p+2}$. Also,
\begin{equation}\label{J^n}
\sup_{n \ge1}
E^W \bigl[ \| J^n_\cdot\|_{L_{p',1}([0,T] \rightarrow V_{p',-\beta
(p,1)})}^{p'}
\bigr]
\le C_T <\infty.
\end{equation}
\end{lemma}

\begin{pf} For any $p \in(1,\infty)$, there exists
$A_p \in(0,\infty)$ such that
\[
\|\mathcal{P}_n v \|_p \le A_p \| v \|_p
     \qquad \mbox{for all $v \in V_{p,0}$}.
\]
(See, e.g., \cite{Gra04}, Theorem 3.5.7, page~213.)
This implies that $\| \mathcal{P}_n v \|_{p,\alpha} \le A_p \| v \|
_{p,\alpha}$
and hence, $\| \mathcal{P}_n v \|_{p',-\alpha} \le A_p \| v \|
_{p',-\alpha}$
for any $p \in(1,\infty)$ and $\alpha\ge0$.
We combine this and (\ref{intb(X)a}) and (\ref{intb(X)b}) to obtain
(\ref{I^n}) and (\ref{J^n}).
\end{pf}

\begin{lemma}\label{Lem.tight1}
For $\beta> \beta(p,1)$,
the laws $\{P^W ( X^n \in\cdot)\}_{n=1}^\infty$
are tight on $C([0,\infty) \rightarrow V_{2 \wedge p',-\beta})$.
\end{lemma}

\begin{pf}
As is easily seen, it is enough to prove the following:\\
(1)
The laws
$\{ P^W ( (X^n_t)_{t \le T} \in\cdot)\}_{n=1}^\infty$
are tight on $C([0,T] \rightarrow V_{2 \wedge p',-\beta})$ for each
fixed $T>0$.
To see (1), we set
\begin{eqnarray*}
\mathcal{S}&=&L_{p,1}\bigl([0,T] \rightarrow V_{p',-\beta(p,1)}\bigr)
+ L_{p',1}([0,T] \rightarrow V_{p',-1})
\\
&&{}+ L_{{2 /\gamma},\gamma}([0,T] \rightarrow V_{2,0}),
  \qquad  \mbox{with }    \gamma\in(0,1/2).
\end{eqnarray*}
We then see from Lemma~\ref{Lem.cpt1} that
\renewcommand{\theequation}{2}
\setcounter{equation}{0}
\begin{eqnarray}
\mathcal{S}\stackrel{\rm compact}{\hookrightarrow} C([0,T]
\rightarrow
V_{2 \wedge p',-\beta}).
\end{eqnarray}
On the other hand, we have that
\renewcommand{\theequation}{\arabic{equation}}
\setcounter{equation}{2}
\begin{eqnarray}
\sup_n E^W \bigl[ \| I^n_\cdot\|_{L_{p,1}([0,T] \rightarrow V_{p',-\beta
(p,1)})}^\delta\bigr]
&\stackrel{\mbox{\scriptsize(\ref{I^n})}}{\le}& C_T <\infty
\qquad\mbox{for some } \delta\in(0,1];\hspace*{-12pt}
\\
\sup_n E^W \bigl[ \| J^n_\cdot\|_{L_{p',1}([0,T] \rightarrow V_{p',-1})} \bigr]
&\stackrel{\mbox{\scriptsize(\ref{J^n})}}{\le}& C_T <\infty;
\\
\quad \sup_n E^W \bigl[ \| X^n_0+W^n_\cdot\|_{L_{2/\gamma
,\gamma}([0,T]
\rightarrow V_{2,0})}\bigr]
&\stackrel{\mbox{\scriptsize(\ref{HBM1/2})}}{\le}& C_T <\infty.
\end{eqnarray}
We conclude from (3)--(5) and (\ref{J+W}) that
\[
\sup_n E^W [ \| X^n_\cdot\|_\mathcal{S}^\delta] \le C_T <\infty
\]
and hence, that for $R>0$,
\setcounter{equation}{5}
\begin{eqnarray}
\sup_n P^W ( \| X^n_\cdot\|_\mathcal{S}> R)
&\le&\frac{1 }{ R^\delta}\sup_n E^W [ \| X^n_\cdot\|_\mathcal
{S}^\delta]\nonumber
\\[-8pt]\\[-8pt]
&\le&\frac{C_T }{ R^\delta} \longrightarrow0   \qquad \mbox{as $R
\longrightarrow\infty$}.\nonumber
\end{eqnarray}
We see from (2) that the set
\[
\{ X_\cdot  ;   \| X^n_\cdot\|_\mathcal{S}\le R \}
\]
is relatively compact in $C([0,T] \rightarrow V_{2 \wedge p',-\beta
})$. Hence, by (6),
we have the tightness (1).
\end{pf}

\begin{lemma}\label{Lem.tight3}
The laws $\{ P^W ( X^n \in\cdot)\}_{n=1}^\infty$
are tight on $L_{\tilde{p}, {\rm loc}}([0,\infty) \rightarrow
V_{\tilde{p},1} )$.
\end{lemma}

\begin{pf} Let $\tilde{p}>1$ and $\tilde{\alpha} >1$ be from
Lemma~\ref{Lem.3.34}.
We may assume that $\tilde{p} \in(1, p'']$. It is enough to prove
the following:
\renewcommand{\theequation}{1}
\setcounter{equation}{0}
\begin{eqnarray}
&&\mbox{The laws
$\{ P^W ( (X^n_t)_{t \le T} \in\cdot)\}_{n=1}^\infty$
are tight on $L_{\tilde{p}}([0,T] \rightarrow V_{\tilde{p},1}
)$}\nonumber
\\[-8pt]\\[-8pt]
&& \mbox{for each fixed
$T>0$.}\nonumber
\end{eqnarray}
To see (1), we set
\[
\mathcal{I}=L_{\tilde{p}}([0,T] \rightarrow V_{\tilde
{p},\tilde{\alpha}})
\cap L_{\tilde{p},\gamma}\bigl([0,T] \rightarrow V_{\tilde
{p},-\beta(p,1)}\bigr)
 \qquad   \mbox{with }    \gamma\in(0,1/2).
\]
Note that
\[
V_{\tilde{p},\tilde{\alpha}} \stackrel{\rm
compact}{\hookrightarrow} V_{\tilde{p},1}
\hookrightarrow V_{\tilde{p},-\beta(p,1)}
\]
and hence, by Lemma~\ref{Lem.cpt2}, that
\renewcommand{\theequation}{2}
\setcounter{equation}{0}
\begin{eqnarray}
\mathcal{I}\stackrel{\rm compact}{\hookrightarrow}
L_{\tilde{p}}([0,T] \rightarrow V_{\tilde{p},1}).
\end{eqnarray}
On the other hand,
\renewcommand{\theequation}{3}
\setcounter{equation}{0}
\begin{eqnarray}
\sup_n E^W \bigl[ \| X^n_\cdot\|_{L_{\tilde{p}}([0,T] \rightarrow
V_{\tilde{p},\tilde
{\alpha}})}\bigr]
\stackrel{\mbox{\scriptsize
(\ref{3.34})}}{\le} C_T <\infty.
\end{eqnarray}
Moreover, for some $\delta\in(0,1]$,
\begin{eqnarray*}
&&\sup_n E^W \bigl[
\| X^n_\cdot\|^\delta_{L_{\tilde{p},\gamma}([0,T] \rightarrow
V_{\tilde{p},-\beta
(p,1)})}\bigr]
\\
&&\hspace*{17pt}\qquad \le \sup_n E^W \bigl[ \| X^n_0+I^n_\cdot+J^n_\cdot
\|_{L_{\tilde{p},\gamma}([0,T] \rightarrow V_{\tilde
{p},-\beta(p,1)})}^\delta\bigr]
\\
&&{}\hspace*{17pt}\qquad\quad+\sup_n E^W \bigl[
\| W^n_\cdot\|^\delta_{L_{\tilde{p},\gamma}([0,T] \rightarrow
V_{2,0})}\bigr]
\\
&&\qquad\stackrel{\mbox{\scriptsize(\ref{HBM1/2})--(\ref{J^n})}}{\le}
C_T <\infty.
\end{eqnarray*}
We conclude from (2) and (3) that
\[
\sup_n E^W [ \| X^n_\cdot\|_\mathcal{I}^\delta] \le C_T <\infty
\]
and hence, that for $R>0$,
\renewcommand{\theequation}{4}
\setcounter{equation}{0}
\begin{eqnarray}
\sup_n P^W ( \| X^n_\cdot\|_\mathcal{I}> R)
&\le&\frac{1 }{ R^\delta}\sup_n E^W [ \| X^n_\cdot\|_\mathcal
{I}^\delta]\nonumber
\\[-8pt]\\[-8pt]
&\le&\frac{C_T }{ R^\delta} \longrightarrow0 \qquad   \mbox{as $R
\longrightarrow\infty$}.\nonumber
\end{eqnarray}
We will see from this and (2) that the set
\[
\{ X_\cdot  ;   \| X^n_\cdot\|_\mathcal{I}\le R \}
\]
is relatively compact in $L_{\tilde{p}}([0,T] \rightarrow
V_{\tilde{p},1})$.
Hence, by (4) we have the tightness~(1).
\end{pf}

Finally, Proposition~\ref{Prop.tight} follows from Lemmas~\ref
{Lem.tight1}, \ref{Lem.tight3} and
the following:

\begin{lemma}\label{Lem.tight0}
Suppose that:

\begin{itemize}[$\bullet$]
\item[$\bullet$]
$(S_j, \rho_j)$ ($j=1,\ldots,m$) are
complete separable
metric spaces such that all of $S_j$ ($j=1,\ldots,m$) are subsets of
a set $S$;

\item[$\bullet$]
$(X_n)_{n \in{\mathbb{N}}}$ is a sequence of random variables
with values in $\bigcap_{j=1}^m S_j$ defined on a probability space
$(\Omega, \mathcal{F}, P)$;

\item[$\bullet$]
$(X_n)_{n \in{\mathbb{N}}}$ is tight in each of
$(S_j, \rho_j)$, $j=1,\ldots,m,$ separately.

Then, there exists a sequence $n(k) \rightarrow\infty$, random variables
$X, \widetilde{X}_k$, $k=1,2,\ldots,$ with values in $\bigcap_{j=1}^m S_j$
defined on a
probability space $(\widetilde{\Omega}, \widetilde{\mathcal{F}},
\widetilde{P})$ such that
\begin{eqnarray*}
\widetilde{P}(\widetilde{X}_k \in\cdot)&=&P\bigl(X_{n(k)}\in\cdot\bigr)
   \qquad \mbox{for all $k=1,2,\ldots;$}
    \\
 \lim_{k \rightarrow\infty}\sum_{j=1}^m\rho_j (X,\widetilde
{X}_k)&=&0  \qquad  \mbox{$\widetilde{P}$-a.s.}
\end{eqnarray*}
\end{itemize}
\end{lemma}

\begin{pf} By induction, it is enough to consider the case of $m=2$.
Let $\varepsilon>0$ be arbitrary. Then, for
$j=1,2$, there exists a compact subset $K_j$ of
$S_j$ such that
\[
P (X_n \in K_j) \ge1 -\varepsilon\qquad
    \mbox{for all $j=1,2$ and $n=1,2,\ldots.$}
\]
Now a very simple but crucial observation is
that $K_1 \cap K_2$ is compact in $S_1 \cap S_2$
with respect to the metric $\rho_1+\rho_2$. Also,
\[
P (X_n \in K_1 \cap K_2) \ge1 -2\varepsilon
   \qquad \mbox{for all $j=1,2$ and $n=1,2,\ldots.$}
\]
These imply that $(X_n)$ is tight in $S_1 \cap S_2$
with respect to the metric $\rho_1+\rho_2$.
Thus, the lemma follows from
Prohorov's theorem (\cite{IW89}, Theorem 2.6, page~7)
and Skorohod's representation theorem
(\cite{IW89}, Theorem 2.7, page~9).
\end{pf}

\begin{remark*} This remark, together with the one after
Lemma~\ref{Lem.conv_intb},
concerns the validity of Theorem~\ref{Thm.SNS} with $p=2$ for all $d$.
Let $\alpha<1$. Then we can also prove that
\renewcommand{\theequation}{\arabic{section}.\arabic{equation}}
\setcounter{equation}{28}
\begin{equation}\label{tight_p,a}
\mbox{the laws
$\{P^W ( X^n \in\cdot)\}_{n=1}^\infty$
are tight on $L_{p'',{\rm loc}}([0,\infty) \rightarrow V_{p,\alpha} )$.}
\end{equation}
This can be seen as follows. We set
\[
\mathcal{I}=L_{p''}([0,T] \rightarrow V_{p,1}) \cap L_{p'',\gamma
}\big([0,T] \rightarrow
V_{p'',-\beta(p,1)}\big),
 \qquad   \mbox{with }    \gamma\in(0,1/2).
\]
Since
\[
V_{p,1} \stackrel{\rm compact}{\hookrightarrow} V_{p,\alpha}
\hookrightarrow V_{p'',-\beta(p,1)},
\]
we have by Lemma~\ref{Lem.cpt2} that
\[
\mathcal{I}\stackrel{\rm compact}{\hookrightarrow}
L_{p''}([0,T] \rightarrow V_{p,\alpha}).
\]
Then we get (\ref{tight_p,a}) by similar argument as in Lemma~\ref
{Lem.tight3}.

By the tightness (\ref{tight_p,a}), Lemmas~\ref{Lem.tight1} and
\ref{Lem.tight0},
we obtain a variant of
Proposition~\ref{Prop.tight} in which
the convergence $\widetilde{X}^k \rightarrow X$, $P$-a.s. takes place
in the metric space
\begin{equation}\label{space_tight2}
C\bigl([0,\infty) \rightarrow V_{2 \wedge p',-\beta} \bigr)
\cap L_{p'',{\rm loc}}\bigl([0,\infty) \rightarrow V_{p,\alpha} \bigr)
\end{equation}
instead of (\ref{space_tight}). We note that this modification of
Proposition~\ref{Prop.tight} is
valid for $p \in[\frac{2d }{ d+2},\infty)$ since we did not use
Lemma~\ref{Lem.3.34}.
\end{remark*}

\section{\texorpdfstring{Proof of Theorems~\protect\ref{Thm.SNS}
and~\protect\ref{Thm.SNS1}}{Proof of Theorems~2.1.3 and~2.2.1}}

\subsection{\texorpdfstring{Proof of Theorem~\protect\ref{Thm.SNS}}{Proof of Theorem~2.1.3}}
\label{verify}

Let $X$ and $\widetilde{X}^k$ be as in Proposition~\ref{Prop.tight}.
We will verify
(\ref{state}) [with $\beta=\beta(p,1)$]
as well as (\ref{X_0=m})--(\ref{SNS}) and (\ref{apriori2}) for $X$.
(\ref{X_0=m}) can easily be seen. In fact,
\begin{eqnarray*}
\widetilde{X}^k_0 & \rightarrow& X_0   \qquad  \mbox{ a.s. in $V_{2
\wedge p',-\beta}$},
\\
\widetilde{X}^k_0 \stackrel{\rm law}{=} X^{n(k)}_0=\mathcal
{P}_{n(k)}\xi& \rightarrow& \xi
   \qquad \mbox{in $V_{2,0}$}.
\end{eqnarray*}
Thus, the laws of $X_0$ and $\xi$ are identical.
\[
\widetilde{X}^k_0\stackrel{\rm law}{=}X^{n(k)}_0=\mathcal
{P}_{n(k)}\xi\rightarrow\xi
   \qquad \mbox{in $V_{2,0}$}.
\]
Note that the function
\[
v_\cdot
  \mapsto
\sup_{t \le T}\| v_t \|_2^2 +\int^T_0\| v_t \|_{p,1}^p\, dt
\]
is lower semi-continuous on the metric space (\ref{space_tight}).
Thus, (\ref{apriori2}) follows from (\ref{apriori}) and
Proposition~\ref{Prop.tight}
via Fatou's lemma.

To show (\ref{X_0indepM}) and (\ref{SNS}), we prepare the following:

\begin{lemma}\label{Lem.conv_intb}
Let $\varphi\in\mathcal{V}$ and $T>0$. Then,
\begin{eqnarray}\label{conv_intb1}
\qquad\hspace*{4pt} \lim_{k \rightarrow\infty}
\int^T_0 |\langle\varphi, (\widetilde{X}^k_t \cdot\nabla)
\widetilde{X}^k_t
-(X_t \cdot\nabla) X_t \rangle|\, dt
&=& 0
  \qquad  \mbox{in probability $(P)$},\hspace*{-6pt}
     \\\label{conv_intb2}
\lim_{k \rightarrow\infty}
\int^T_0 |\langle e( \varphi), \tau(\widetilde{X}^k_t)-\tau( X_t)
\rangle|\,dt
&=&0
 \qquad   \mbox{in $L_1(P)$},
    \\\label{conv_intb}
\lim_{k \rightarrow\infty}
\int^T_0 \big\langle\varphi, \mathcal{P}_{n(k)}b (\widetilde{X}^k_t)
-b ( X_t) \big\rangle \,dt
&=& 0
  \qquad  \mbox{in probability $(P)$}.
\end{eqnarray}
\end{lemma}

\begin{pf} We write $Z^k_t =\widetilde{X}^k_t -X_t$ to simplify the notation.
We start by proving that
\begin{equation}\label{convp_1}
\lim_{k \rightarrow\infty}E \biggl[
\int^T_0 \| Z^k_t \|_{p_1,1}^{p_1} \,dt \biggr]=0,   \qquad  \mbox{if $p_1 <p$}.
\end{equation}
By Proposition~\ref{Prop.tight},
\[
I_k \stackrel{\rm def.}{=}\int^T_0 \| Z^k_t \|_{1,1} \,dt \stackrel{k
\rightarrow\infty}\longrightarrow0,
    \qquad  \mbox{$P$-a.s.}
\]
Moreover, the random variables $\{I_k\}_{k \ge1}$ are uniformly
integrable since
\[
E [ I_k^p]
\stackrel{\mbox{\scriptsize(\ref{apriori})}}{\le} C_T <\infty.
\]
Therefore,
\renewcommand{\theequation}{2}
\setcounter{equation}{0}
\begin{eqnarray}
\lim_{k \rightarrow\infty}E [ I_k ]=0.
\end{eqnarray}
Let $k (m) \nearrow\infty$ be such that
\renewcommand{\theequation}{3}
\setcounter{equation}{0}
\begin{eqnarray}
\Phi_{m,t}\stackrel{\rm def.}{=}\bigl|Z^{k(m)}_t \bigr|+\bigl|\nabla Z^{k(m)}_t
\bigr|\stackrel{m \rightarrow\infty
}{\longrightarrow}0,\qquad
\mbox{$dt|_{[0,T]} \times dx \times P$-a.e.},
\end{eqnarray}
where $dt|_{[0,T]} \times dx$ denotes the Lebesgue measure on $[0,T]
\times{\mathbb{T}}^d$. Such a sequence
$k (m)$ exists by (2). The sequence $\{ \Phi_{m, \cdot} \}_{m \ge1}$
is uniformly
integrable with respect to $dt|_{[0,T]} \times dx \times P$. In fact,
\[
E \biggl[ \int^T_0  \hspace*{-2pt}\int_{{\mathbb{T}}^d}\Phi_{m, t}^p \,dt \biggr]
\stackrel{\mbox
{\scriptsize(\ref{apriori})}}{\le} C_T <\infty.
\]
Therefore, (3), together with this uniform integrability, implies (\ref
{convp_1}) along
the subsequence $ k (m)$.
Finally, we get rid of the subsequence, since the subsequence as $k(m)$
above can be
chosen from any subsequence of $k$ given in advance.
We now prove (\ref{conv_intb1}).
Since
\[
(\widetilde{X}^k_t \cdot\nabla) \widetilde{X}^k_t -(X_t \cdot
\nabla)X_t
=(Z^k_t\cdot\nabla)\widetilde{X}^k_t
+(X_t \cdot\nabla) Z^k_t,
\]
we have
\[
\int^T_0 |\langle\varphi, (\widetilde{X}^k_t \cdot\nabla
)\widetilde{X}^k_t
-(X_t \cdot\nabla) X_t \rangle|\,dt \le J_1+J_2,
\]
where
\[
J_1 =
\int^T_0 |\langle\varphi, (Z^k_t \cdot\nabla)\widetilde{X}^k_t
\rangle|\,dt\quad     \mbox{and}\quad
J_2
= \int^T_0 |\langle\varphi, (X_t \cdot\nabla)Z^k_t \rangle|\,dt.
\]
We may take $p_1$ in (\ref{convp_1}) as bigger than $\frac{3d }{ d+2}$
so that
there exists $0<\alpha<1$ such that
$\frac{2d }{ d+2 \alpha} <p_1$.
Then by (\ref{Buvw<a2}), we have that
\[
|\langle\varphi, ( Z^k_t \cdot\nabla)\widetilde{X}^k_t \rangle|
\le C\| Z^k_t \|_{p_1,\alpha}\| \widetilde{X}^k_t \|_2 \|\varphi\|
_{p_1,\beta(p_1,\alpha)}
\]
and hence that
\[
J_1 \le
C\|\varphi\|_{p_1,\beta(p_1,\alpha)}\sup_{t \le T}\| \widetilde
{X}^k_t \|_2
\int^T_0 \| Z^k_t \|_{p_1,\alpha}\,dt.
\]
By (\ref{apriori}) and (\ref{convp_1}),
\[
\sup_{k \ge1}E\Bigl[\sup_{t \le T}\| \widetilde{X}^k_t \|_2^2\Bigr]<\infty
  \quad  \mbox{and}\quad
\lim_{k \rightarrow\infty}\int^T_0 \| Z^k_t \|_{p_1,\alpha}\,dt=0
  \qquad  \mbox{$P$-a.s.}
\]
Thus, $\lim_{k \rightarrow\infty}J_1 =0$ in probability.
On the other hand, we have by (\ref{Buvw<a3}) that
\[
|\langle\varphi, (X_t \cdot\nabla) Z^k_t \rangle|
\le C\| Z^k_t \|_{p_1,\alpha}\| X_t \|_2 \|\varphi\|_{p_1,\beta
(p_1,\alpha)}
\]
and hence that
\[
J_2 \le
C\|\varphi\|_{p_1,\beta(p_1,\alpha)}\sup_{t \le T}\| X_t \|_2
\int^T_0 \| Z^k_t \|_{p_1,\alpha}\,dt.
\]
By (\ref{apriori2}) and (\ref{convp_1}),
\[
E\Bigl[\sup_{t \le T}\| X_t \|_2^2\Bigr]<\infty
   \quad \mbox{and}\quad
\lim_{k \rightarrow\infty}\int^T_0 \| Z^k_t \|_{p_1,\alpha}\,dt=0
  \qquad  \mbox{$P$-a.s.}
\]
Thus, $\lim_{k \rightarrow\infty}J_2 =0$ in probability.

We now turn to (\ref{conv_intb2}). It is enough to prove that
\renewcommand{\theequation}{4}
\setcounter{equation}{0}
\begin{eqnarray}
\lim_{k \rightarrow\infty}E \biggl[
\int^T_0 \| \tau(\widetilde{X}^k_t) - \tau(X_t) \|_1 \,dt
\biggr]=0 .
\end{eqnarray}
Again, let $k (m)$ be such that (3) holds. Then,
\renewcommand{\theequation}{5}
\setcounter{equation}{0}
\begin{eqnarray}
\lim_{m \rightarrow\infty}\tau\bigl( \widetilde{X}^{k(m)}_t\bigr) =\tau(
X_t ),\qquad
\mbox{$dt|_{[0,T]} \times dx \times P$-a.e.}
\end{eqnarray}
On the other hand, we have for $p'=\frac{p }{ p-1}$ that
\[
E \biggl[ \int^T_0 dt \int_{{\mathbb{T}}^d}|\tau(\widetilde
{X}^k_t)|^{p'} \biggr]
\le CE \biggl[ \int^T_0 dt \int_{{\mathbb{T}}^d} \bigl( 1+|e
(\widetilde{X}^k_t)|
\bigr)^p \biggr]
\stackrel{\mbox{\scriptsize(\ref{apriori})}}{\le} C_T <\infty,
\]
which implies that $\tau(\widetilde{X}^k_t)$, $k \in{\mathbb{N}}$
are uniformly integrable with respect to $dt|_{[0,T]} \times dx \times P$.
Therefore, (5), together with this uniform integrability, implies (4)
along the subsequence $k(m)$.
Finally, we get rid of the subsequence, since the subsequence as $k
(m)$ above can be
chosen from any subsequence of $k$ given in advance.

Equation~(\ref{conv_intb}) follows from (\ref{conv_intb1}) and
(\ref{conv_intb2}). Since $\varphi\in\mathcal{V}$ is fixed and $k$
is tending to
$\infty$, we do not have to care about $\mathcal{P}_{n(k)}$ here.
\end{pf}

\begin{remark*} If $p=2$, then Lemma~\ref{Lem.conv_intb} is
valid for
all $d$. This is for the following reason.
By inspection of the proof above, we see immediately that (\ref
{conv_intb1}) follows also from
the modification of Proposition~\ref{Prop.tight} mentioned at the end
of Section \ref
{conv_appro}. Also, for $p=2$, (\ref{conv_intb2})
is equivalent to
\[
\lim_{k \rightarrow\infty}
\int^T_0 \langle\Delta\varphi, \widetilde{X}^k_t- X_t \rangle\, dt=0
  \qquad  \mbox{in $L_1(P)$},
\]
which also follows from the modification of Proposition~\ref
{Prop.tight} mentioned at
the end of Section \ref{conv_appro}.
\end{remark*}

\begin{lemma}\label{Lem.Y=BM}
Let
\renewcommand{\theequation}{\arabic{section}.\arabic{equation}}
\setcounter{equation}{4}
\begin{equation}\label{Y_t}
Y_t =Y_t (X)=X_t-X_0-\int^t_0 b(X_s)\,ds, \qquad    t \ge0.
\end{equation}
Then, $Y_\cdot$ is a ${\rm BM}(V_{2,0}, \Gamma)$.
Moreover, $Y_{t+\cdot}-Y_t$ and
$\{ \langle\varphi, X_s \rangle  ;   s \le t, \varphi\in\mathcal
{V}\}$
are independent for any $t \ge0$.
\end{lemma}

\begin{pf}
It is enough to prove that for each
$\varphi\in\mathcal{V}$ and $0 \le s <t$,
\renewcommand{\theequation}{1}
\setcounter{equation}{0}
\begin{eqnarray}
E [ \exp(\mathbf{i}\langle\varphi, Y_t -Y_s\rangle
)| \mathcal{G}_s ]
=\exp\biggl(-\frac{t-s }{2}\langle\varphi, \Gamma\varphi\rangle
\biggr),
   \qquad \mbox{a.s.},
\end{eqnarray}
where $\mathcal{G}_s=\sigma(\langle\varphi, X_u \rangle  ;   u
\le s, \varphi\in\mathcal{V})$.
We set
\[
F (X)=f (\langle\varphi_1 ,X_{u_1}\rangle, \ldots,\langle\varphi_n
,X_{u_n}\rangle),
\]
where $f \in C_{\rm b}({\mathbb{R}}^n)$,
$0 \le u_1<\cdots<u_n \le s$ and $\varphi_1,\ldots,\varphi_n \in\mathcal
{V}$ are chosen
arbitrarily in advance. Then (1) can be verified by showing that
\renewcommand{\theequation}{2}
\setcounter{equation}{0}
\begin{eqnarray}
E [ \exp(\mathbf{i}\langle\varphi, Y_t -Y_s\rangle)
F(X) ]
=\exp\biggl(-\frac{t-s }{2}\langle\varphi, \Gamma\varphi\rangle
\biggr) E[F(X)].
\end{eqnarray}
Let
\[
Y^k_t =\widetilde{X}^k_t-\widetilde{X}^k_0-\int^t_0
\mathcal{P}_{n(k)}b(\widetilde{X}^k_s)\,ds,  \qquad   t \ge0.
\]
Then we see from Theorem~\ref{Thm.Gar} that
\renewcommand{\theequation}{3}
\setcounter{equation}{0}
\begin{eqnarray}
&&E [ \exp(\mathbf{i}\langle\varphi, Y^k_t -Y^k_s\rangle
) F(\widetilde
{X}^k) ]\nonumber
\\[-8pt]\\[-8pt]
&&\qquad =\exp\biggl(-\frac{t-s }{2}\bigl\langle\varphi, \Gamma\mathcal{P}_{n(k)}
\varphi\bigr\rangle\biggr)
E[F(\widetilde{X}^k)].\nonumber
\end{eqnarray}
Moreover, we have
\[
\lim_{k \rightarrow\infty}\langle\varphi, Y^k_t -Y^k_s\rangle
\stackrel{\mbox{\scriptsize(\ref{tlX->X}),(\ref{conv_intb})}}{=}
\lim_{k \rightarrow\infty}\langle\varphi, Y_t -Y_s\rangle
    \qquad  \mbox{in probability}
\]
and hence,
\[
\lim_{k \rightarrow\infty}\mbox{LHS of (3)}=\mbox{ LHS of (2)}.
\]
On the other hand,
\[
\lim_{k \rightarrow\infty}\mbox{ RHS of (3)}
\stackrel{\mbox{\scriptsize(\ref{tlX->X})}}{=}
\mbox{ RHS of (2)}.
\]
These prove (2).
\end{pf}

Finally, we prove (\ref{state}) with $\beta=\beta(p,1)$.
It follows from (\ref{apriori2}) that
\[
X \in L_{p,{\rm loc}} \bigl([0,\infty) \rightarrow V_{p,1}\bigr)
\cap L_{\infty,{\rm loc}} \bigl([0,\infty) \rightarrow V_{2,0}\bigr).
\]
Thus, it remains
to show that $X \in C([0,\infty) \rightarrow V_{2 \wedge p',-\beta(p,1)})$.
But this follows from
Lemma~\ref{Lem.intb(X)} and that $Y \in C([0,\infty) \rightarrow V_{2,0})$.

\subsection{\texorpdfstring{Proof of Theorem~\protect\ref{Thm.SNS1}}{Proof of Theorem~2.2.1}}

Here we can follow the argument of \cite{MNRR96}, Theorem 4.29, page~254,
almost verbatim. We will present it for the convenience of the readers.

We need two technical lemmas.

\begin{lemma}\label{Lem.1.2}
Let $H$ be a Hilbert space and $V$ be a Banach space such that
\[
V \hookrightarrow H \hookrightarrow V^*.
\]
Suppose that $f \in L_p ([0,T] \rightarrow V)$ ($p \in(1,\infty)$,
$T>0$) has
derivative $f'$ in $L_{p'} ([0,T] \rightarrow V^*)$. Then,
\renewcommand{\theequation}{\arabic{section}.\arabic{equation}}
\setcounter{equation}{5}
\begin{equation}\label{1.2}
\frac{d }{ dt}|f|_H^2=2{}_V\langle f, f' \rangle_{V^*}
\end{equation}
in the distributional sense on $(0,T)$.
\end{lemma}

\begin{pf} The case of $p=2$ can be found in \cite{Te79}, Lemma 1.2,
pages~60--61. The extension to general $p$ is straightforward.
\end{pf}

\begin{lemma}[(\cite{MNRR96}, Lemma 4.35, page~255)]\label{Lem.4.35}
Let $q \in(2,\infty)$ if $d=2$ and $q \in[2,\frac{2d }{ d-2}]$ if $d
\ge3$.
Then there exists $c \in(0,\infty)$ such that
\setcounter{equation}{6}
\begin{equation}\label{4.35}
\| v\|_q \le c\| v\|_2^\theta\| \nabla v\|_2^{1-\theta}
\qquad \mbox{with } \theta=\frac{2d-q(d-2) }{2q}
\end{equation}
for all $v \in V_{2,1}$ with $\int_{{\mathbb{T}}^d}v =0$.
\end{lemma}

Let $X$ and $\widetilde{X}$ be as in the assumptions of Theorem~\ref
{Thm.SNS1} and
\[
Z_t=X_t-\widetilde{X}_t=\int^t_0 \bigl( b (X_s)-b (\widetilde
{X}_s)\bigr)\,ds.
\]
Then,
%
\renewcommand{\theequation}{1}
\setcounter{equation}{0}
\begin{eqnarray}
Z_\cdot\in L_{p, {\rm loc}}\bigl([0,\infty) \rightarrow V_{p,1}\bigr)
\end{eqnarray}
%
and by Lemma~\ref{Lem.intb(X)},
\renewcommand{\theequation}{2}
\setcounter{equation}{0}
\begin{eqnarray}
\partial_tZ_\cdot= b (X_\cdot)-b (\widetilde{X}_\cdot) \in
L_{p, {\rm loc}}\bigl([0,\infty) \rightarrow V_{p',-\beta(p,1)}\bigr).
\end{eqnarray}
 Since $p \ge p'$ and $\beta(p,1)=1$
for $p \ge1+\frac{d }{2} (\ge\frac{4d }{ d+2})$, we see from (2) and
Lemma~\ref{Lem.1.2} (applied to $f=Z_\cdot$ and $V=V_{p,1}$)
that
%
\renewcommand{\theequation}{3}
\setcounter{equation}{0}
\begin{eqnarray}
\frac{1}{2}\frac{d }{ dt}\| Z_t\|_2^2
\stackrel{\mbox{\scriptsize(\ref{1.2})}}{=}
\langle Z_t, b (X_t )-b (\widetilde{X}_t ) \rangle=-I_t-J_t
\end{eqnarray}
in the distributional sense, where
\[
I_t = \langle Z_t, (X_t \cdot\nabla)X_t- (\widetilde{X}_t \cdot
\nabla)\widetilde{X}_t
\rangle
  \quad  \mbox{and}\quad
J_t = \langle e (Z_t), \tau(X_t )-\tau(\widetilde{X}_t ) \rangle.
\]
We have by \cite{MNRR96}, formula (1.25), page~198 and formula (1.11), page~196,
that
%
\renewcommand{\theequation}{4}
\setcounter{equation}{0}
\begin{eqnarray}
J_t \ge c_1 \| e (Z_t)\|_2^2 \ge c_2 \| \nabla Z_t \|_2^2.
\end{eqnarray}
On the other hand, since $\widetilde{X}_t=X_t-Z_t$, we see that
\[
\langle Z_t, (\widetilde{X}_t \cdot\nabla)\widetilde{X}_t \rangle
\stackrel{\mbox{\scriptsize(\ref{uu=0})}}{=}
\langle Z_t, (\widetilde{X}_t \cdot\nabla)X_t \rangle
=\bigl\langle Z_t, \bigl((X_t-Z_t) \cdot\nabla\bigr)X_t \bigr\rangle,
\]
and hence that
\[
I_t=\langle Z_t, (Z_t \cdot\nabla)X_t \rangle.
\]
Therefore,
%
\renewcommand{\theequation}{5}
\setcounter{equation}{0}
\begin{eqnarray}
|I_t| & \stackrel{{1 / p}+{(p-1)/(2p)}+{(p-1) /(2p)}=1}{\le}&
\| \nabla X_t\|_p\|Z_t\|_{{2p / (p-1)}}^2\nonumber
\\
&\stackrel{\mbox{\scriptsize(\ref{4.35})}}{\le}&\hspace*{-50pt}
C_3 \| \nabla Z_t\|_2^{{d / p}}\| \nabla X_t\|_p\| Z_t\|_2^{{(2p-d)
/ p}}
\\
& \stackrel{{{d /(2p)}}+{(2p-d) /(2p)}=1}{\le} &\hspace*{-15pt}
c_2\| \nabla Z_t\|_2^2 +C_4\| \nabla X_t\|_p^{{2p /(2p-d)}} \| Z_t\|_2^2.\nonumber
\end{eqnarray}
We see from (3)--(5) that
\[
\frac{1}{2}\frac{d }{ dt}\| Z_t\|_2^2
\le C_4\| \nabla X_t\|_p^{{2p /(2p-d)}} \| Z_t\|_2^2.
\]
Since $\frac{2p }{2p-d} \le p$, this implies via Gronwall's lemma
(we need an appropriate generalization since the derivative above is
in the distributional sense)
that
\[
\| Z_t\|_2^2 \le\| Z_0\|_2^2\exp
\biggl( C_4\int^t_0 \| \nabla X_s\|_p^{{2p /(2p-d)}}\,ds \biggr).
\]
This proves that $\| Z_t\|_2\equiv0$.

\section*{Acknowledgments}

The authors thank Professor Josef M\'alek for valuable
comments on the earlier version of this article.
N. Yoshida thanks Professor Kenji Nakanishi for useful conversations.

%

\printaddresses

\end{document}